\documentclass{amsart}
\usepackage{amsmath, amssymb, amsthm}
\usepackage[utf8]{inputenc}
\usepackage[T1]{fontenc}
\usepackage{graphicx}
\usepackage[colorlinks=true, allcolors=blue]{hyperref}
\usepackage{comment}
\usepackage{amsmath}
\usepackage{amssymb}
\IfFileExists{stmaryrd.sty}{\usepackage{stmaryrd}}{\newcommand{\shortrightarrow}{\rightarrow}}
\newcommand{\calL}{\mathcal{L}}
\newtheorem{theorem}{Theorem}[section]
\newtheorem{lemma}[theorem]{Lemma}
\newtheorem{definition}[theorem]{Definition}
\newtheorem{corollary}[theorem]{Corollary}
\newtheorem{example}[theorem]{Example}
\newtheorem{proposition}[theorem]{Proposition}

\newcommand{\lan}{\langle}
\newcommand{\ran}{\rangle}
\usepackage{xcolor}

\usepackage{adjustbox}

\begin{document}
	\title[Three-dimensional Jordan algebras]{Isomorphism classes and geometric classification of three-dimensional Jordan algebras and the length function}
	
	\author[Hern\'andez]{Isabel Hern\'andez}
	\address{CONAHCYT-CIMAT Unidad M\'erida}
	\email{isabel@cimat.mx}

	\author[Martin]{ Mar\'ia Eugenia Martin}
	\address{Universidade Federal do Paran\'a}
	\email{eugeniamartin@gmail.com}
	
	\author[Rodrigues]{Rodrigo Lucas Rodrigues}
	\address{Department of Mathematics\\ Federal University of Cear\'a\\ Pici Campus, Block 914\\ 60455-760, Fortaleza, CE Brazil}
	\email{rodrigo@mat.ufc.br}
	
	\author[Silva]{Carlos Andr\'e Gomes Silva}
	\address{Instituto de Matem\'atica e Estat\'istica \\ Universidade de S\~ao Paulo, \\ S\~ao Paulo - SP, Brazil} 
	\email{carlosmat@ime.usp.br}
	\begin{abstract}
		We investigate the algebraic classification of Jordan algebras of dimensions at most three over arbitrary fields with more than two elements. Our approach combines the theory of nilalgebras and radicals with Peirce, Wedderburn decompositions and the length function. In characteristic two, we adapt the Peirce decomposition to organize the classification, including the field-dependent isomorphism conditions among parametrized families. We also study the varieties of Jordan algebra structures up to dimension three over an algebraically closed field of characteristic two. Besides explicit degenerations and standard arguments, we introduce the alicyclic and characteristic conditions as geometric obstructions to degeneration. We determine the rigid algebras and irreducible components in both dimensions. In dimension three, one component is the nine-dimensional affine space of square-zero products; unlike the other three components, it contains no open orbit.
	\end{abstract}
	\maketitle
	\section{Introduction}
	Jordan algebras originated in 1932 in the study of the mathematical foundations of quantum mechanics, when Pascual Jordan introduced a class of commutative nonassociative algebras \cite{Jordan1932} in an attempt to formalize aspects of quantum mechanics. This approach was further developed by Jordan, von Neumann and Wigner in 1934, classifying finite-dimensional formally real Jordan algebras \cite{JordanVonNeumannWigner1934}. Expanding these works, Albert obtained fundamental results on their structure, including Peirce decompositions, the study of radicals and a version of the Wedderburn principal theorem \cite{Albert1946,AlbertWedderburn1947,Albert1947}. These investigations were followed by the classification and representation of semisimple Jordan algebras by F. Jacobson and N. Jacobson \cite{JacobsonJacobson1949}. Later, in 1983, Zelmanov worked on prime Jordan algebras providing a substantial advance in the structural classification of Jordan algebras beyond the finite-dimensional semisimple condition \cite{Zelmanov1983}.\\
	
	The algebraic classification of low-dimensional Jordan algebras has been historically developed in parallel to the classification of other classes of nonassociative algebras. In 1975, Gerstenhaber and Myung classified four-dimensional commutative power-associative nilalgebras over fields of characteristic different from two \cite{GerstenhaberMyung1975}. Another significant advance was made by Petersson in 2000, who classified all two-dimensional nonassociative algebras over arbitrary fields \cite{Petersson2000}. Nevertheless, extending classifications to higher dimensions introduces additional difficulties, particularly concerning field-dependent restrictions such as the existence of roots for quadratic and cubic polynomials. Kashuba and Shestakov investigated three-dimensional Jordan algebras over algebraically closed fields placing restrictions on the characteristic of the field \cite{KashubaShestakov} and the arithmetic reasoning posed an obstacle to extend the presented method to arbitrary fields of positive characteristic. In 2016, Hegazi and Abdelwahab classified nilpotent Jordan algebras of dimension at most three over arbitrary fields \cite{14}; however, this did not cover nonnilpotent Jordan algebras completely, as it could not include nilalgebras that need not be nilpotent in characteristic two. More recently, Cedilnik and Jerman classified three-dimensional zeropotent algebras over arbitrary fields \cite{19}, covering an important subclass of Jordan algebras in characteristic two, but leaving algebras with nonzero squares beyond its scope. Recently, Diniz, Gonçalves, da Silva and Souza obtained a classification of two-dimensional Jordan algebras over arbitrary fields of characteristic different from two \cite{8}. These developments reveal the remaining difficulties in obtaining a unified low-dimensional classification over general fields, especially in characteristic two, where the usual linearization of polynomial identities may not provide a complete linearization, classical structural decompositions require reconsideration, and isomorphism classes can depend essentially on the arithmetic of the ground field which need not be algebraically closed.\\
	
	The theory of deformations of algebras was initiated by Gerstenhaber in 1964, who developed a framework for studying deformations and rigidity of associative algebras \cite{Gerstenhaber1964}. In his work, Flanigan consolidated the geometric aspects of the theory, investigating varieties of algebra structures and orbit closures \cite{Flanigan1968}. In 2005, Alberca Bjerregaard, Loos and Martín González investigated derivations and automorphisms of Jordan algebras in characteristic two \cite{AlbercaLoosMartin2005}, although the degeneration relations and irreducible components of the low-dimensional varieties remained outside their investigation their work pioneered the investigation of the geometric aspects of the variety of Jordan algebras in characteristic two. Subsequently, Kashuba studied the varieties of unital Jordan algebras of dimension four and five in 2006, determining their irreducible components \cite{Kashuba2006}, but without encompassing nonunital algebras or characteristic two. In 2007, Kashuba and Shestakov investigated the geometric classification of three-dimensional Jordan algebras \cite{KashubaShestakov}, leaving characteristic two out of the scope of their work. In 2018, Kashuba and Martin classified the variety of five-dimensional nilpotent Jordan algebras, showing that one of its five irreducible components contains no dense orbit and is instead determined by an infinite family of nonisomorphic algebras \cite{KashubaMartin2018}. However, their investigation was restricted to nilpotent algebras over fields of characteristic different from two and three. \\
	
	The geometric analysis of the present work therefore develops additional tools instead of relying only on previously established degeneration arguments. For instance, an algebra of dimension $n$ is called \emph{alicyclic} when, for every element $x$, the vectors $x,x^2,\ldots,x^n$ are linearly dependent. We show that this is a Zariski-closed, isomorphism-invariant condition and hence is preserved under degeneration. In dimension three we also consider polynomial relations between coefficients of the characteristic polynomial of the multiplication operator $L_x$, the resulting \emph{characteristic} conditions are likewise preserved under degeneration. Together with polynomial identities, annihilator and automorphism dimensions, and explicit one-parameter changes of basis, these conditions supply obstructions to nonexistent degenerations. \\
	
	The paper is organized as follows. In Section 3 we classify commutative associative algebras of dimensions one, two and three. The treatment separates nilalgebras from algebras admitting idempotents, this is essentially the same reasoning that Goodaire and Zhou \cite{GoodaireZhou2000} have used for the study of alternative rings, uses Peirce decomposition and the Wedderburn-Artin theory where appropriate, and gives isomorphism conditions for the parametrized families. Section 4 studies Jordan algebras, isolating the cases that are not associative and treating characteristic two separately, including the contribution of zeropotent structures in dimension three. In Section 5 we use Peirce decompositions to complete the analysis of three-dimensional unital and nonunital Jordan algebras. The length and the introduced \emph{polynomial degree} provide complementary ways of sorting associative cases, the resulting classification is again followed by criteria for identifying isomorphic parametrized algebras.\\
	
	Finally, Section 6 concerns the geometric classification over algebraically closed fields of characteristic two. We establish the degeneration criteria described above and determine the orbit closures in dimension two. In dimension three, we determine the irreducible components and the rigid algebras. The resulting picture highlights both the effectiveness of geometric obstructions to degeneration and the ways in which Jordan algebra varieties in characteristic two differ from the more familiar settings.

	\section{Preliminaries}
	Throughout this paper, $K$ denotes a field with at least three elements. An algebra $A$ is called simple whenever $A^2 \neq 0$ and the only ideals of $A$ are $\{0\}$ and $A$. \\
	
	An algebra is called Artinian (resp. left)  if it satisfies the descending chain condition for ideals (resp. left) , that is, if $I_1 \supseteq I_2 \supseteq \cdots$ is such a chain, then there exists $k \geq 1$ such that $I_n = I_k$, for every $n \geq k$. It is clear that every finite-dimensional algebra is Artinian.  \\
	
	An algebra $A$ is said to be power-associative if any subalgebra generated by an element of $A$ is associative. This is equivalent to defining the powers of $x \in A$ recursively by $x^1 = x$, $x^{i+1} = xx^{i}$, and requiring that $x^i x^j = x^{i+j}$, for all $x \in A$, for every pair of positive integers $(i,j)$.\\
	
	An element $x$ of a power-associative algebra $A$ is called nilpotent if there exists an integer $m > 1$ such that $x^m=0$. An algebra (ideal) $A$ consisting only of nilpotent elements is called a nilalgebra (nilideal). We say that a nilalgebra $A$ has nilindex $n \geq 2$, and denote it by ${\rm nil}(A) = n$, if $x^n = 0$ for all $x \in A$ and there exists an element $y \in A$ such that $y^{n-1} \neq 0$. We define $\text{nil}(0)=1$, where $0$ denotes the algebra $A=\{0\}$. \\
	
	An algebra (or an ideal) is called nilpotent if there exists an integer $m \geq 1$ 
	such that every product of $m$ elements with any distribution of parentheses vanishes.
	
	\begin{proposition}\label{nilalgebra}\cite[Theorem 3.2]{11}
		Any finite-dimensional alternative (in particular, associative) nilalgebra is nilpotent.    
	\end{proposition}
	
	\begin{proposition}\cite[§4 III, IV]{15}\label{propFrobenius}
		Let $A$ be an associative nilpotent algebra. If $xy=x$ for some $x,y\in A$, then $x=0$. Consequently, if $x_1,\ldots,x_n\in A$ and $x_1x_2\ldots x_n\neq 0$ then the set $\{x_1,x_1x_2,\ldots,x_1x_2\ldots x_n\}$ is linearly independent.
	\end{proposition}
	
	\begin{corollary}\label{nilbasis3}
		Let $A\neq 0$ be a finite-dimensional power-associative nilalgebra. Then ${\rm nil}(A) \leq {\rm dim}(A) + 1$.  
	\end{corollary}
	\begin{proof}
		If some element $a\in A$ satisfies $a^{{\rm dim}(A)+1}\neq 0$, by Proposition \ref{nilalgebra}, we may apply Proposition \ref{propFrobenius} to the subalgebra generated by $a$ and conclude that the set $\{a,a^2,\ldots,a^{{\rm dim}(A)+1}\}$ is linearly independent, a contradiction.
	\end{proof}
	
	\begin{lemma}\label{nilbasis2}
		For every power-associative nilalgebra $A$ of nilindex $2\leq m\leq {\rm dim}(A)+1$, there is an element $x \in A$ such that $\{x, x^2, \ldots, x^{m-1}  \}$ is a linearly independent set of $A$.
	\end{lemma}
	\begin{proof}
		This is a direct application of Proposition \ref{propFrobenius} to the associative algebra generated by an element $a\in A$ satisfying $a^m=0$ and $a^{m-1}\neq 0$.
	\end{proof}
	
	Let $A$ be a finite-dimensional power-associative algebra over a field $K$. It is known that $A$ has a unique maximal nilideal ${\rm Rad}(A)$ such that the quotient algebra $A/{\rm Rad}(A)$ has maximal nilideal 0. The ideal ${\rm Rad}(A)$ is called the nilradical of $A$.
	
	\begin{theorem}\cite{16}
		Let $A$ be a finite-dimensional associative algebra. Then
		\begin{itemize}
			\item[(i)] The nilradical ${\rm Rad}(A)$ is nilpotent;
			\item[(ii)] The quotient algebra $A/{\rm Rad}(A)$ is isomorphic to a finite direct sum of simple associative algebras.
		\end{itemize}
	\end{theorem}
	The Jacobson radical of a unital finite-dimensional associative algebra $A$, denoted by $J(A)$, is the intersection of all the left maximal ideals of $A$. If $A$ is an Artinian algebra, $J(A)$ is nilpotent and it is the sum of all left nilpotent ideals of $A$. In particular, ${\rm Rad}(A) \subset J(A)$. On the other hand, since $A/{\rm Rad}(A)$ does not contain nonzero nilideals and $J(A)/{\rm Rad}(A)$ is a nilideal of $A/{\rm Rad}(A)$, it follows that $J(A) \subset {\rm Rad}(A)$. Therefore, $J(A) = {\rm Rad}(A)$.

	\begin{theorem}\cite{17} (Wedderburn-Artin) 
		Let $A$ be a unital finite-dimensional associative algebra. Then there exist positive integers $n_1, \ldots, n_r$ and division algebras $D_1, \ldots, D_r$ such that $A/J(A) \cong M_{n_1}(D_1) \oplus \cdots \oplus M_{n_r}(D_r)$, where $M_{n_i}(D_i)$ is a matrix algebra with coefficients in $D_i$, for every $1 \leq i \leq r$. Furthermore, if $A/J(A) \cong M_{t_1}(D_1') \oplus \cdots \oplus M_{t_s}(D_s')$ is another decomposition, then $r = s$, $n_i = t_i$ and $D_i \cong D_i'$ ($i = 1, \ldots, r$, reordering the indexes if necessary.)   
	\end{theorem}

	\begin{proposition}\cite[Theorem 4.3]{11}\label{NilJordan}
		Any finite-dimensional Jordan nilalgebra over a field of characteristic different from two is nilpotent.
	\end{proposition}
	
	An element $e$ of an (arbitrary) algebra $A$ is called an idempotent in the case $e^2 = e \neq 0$.\begin{proposition}\label{idempotent}\cite[Proposition 3.3]{11}
		Any finite-dimensional power-associative algebra, which is not a nilalgebra, contains an idempotent.
	\end{proposition}
	\begin{lemma}\label{newpeirce}
		Any commutative associative algebra $A$ with an idempotent $e$ has the Peirce decomposition $A=A_1 (e) \oplus A_0 (e)$, where $A_1 (e) = \{x \in A \ | \ ex = x\}$ and $A_0 (e) = \{x \in A \ | \ ex = 0\}$ are orthogonal subalgebras of $A$. In this case, we say that $A$ has type $({\rm dim}(A_1(e)), {\rm dim}(A_0(e)))$. 
	\end{lemma}
	\begin{proof}
		It suffices to decompose an arbitrary element $x\in A$ as $x = ex + (x - ex)$.
	\end{proof}
	\begin{theorem}\cite{11}[Peirce decomposition for Jordan algebras over a field of characteristic different from two]
		Let $J$ be a Jordan algebra over a field $K$ with $\text{char}(K)\neq 2$, and let $e\in J$ be an idempotent. Then $J$ decomposes as a direct sum
		$$
		J = J_1(e)\;\oplus\; J_{1/2}(e)\;\oplus\; J_0(e),
		$$
		where
		$$
		J_\lambda(e)=\{x\in J \mid ex=\lambda x\}, \qquad 
		\lambda\in\left\{1,\tfrac12,0\right\}.
		$$
		Moreover, the Peirce components satisfy the multiplication rules
		$$
		\begin{aligned}
			&J_1(e)\,J_1(e)\subseteq J_1(e), \qquad
			J_0(e)\,J_0(e)\subseteq J_0(e), \qquad
			J_1(e)\,J_0(e)=0,\\
			&J_1(e)\,J_{1/2}(e)\subseteq J_{1/2}(e), \qquad
			J_0(e)\,J_{1/2}(e)\subseteq J_{1/2}(e), \qquad
			J_{1/2}(e)\,J_{1/2}(e)\subseteq J_1(e)\oplus J_0(e).
		\end{aligned}
		$$
	\end{theorem}
	Over a field with at least three elements, any Jordan algebra is power associative \cite{1,4} so that over a field of characteristic 2 with at least four elements a linearization of $a^2(ba)+(a^2b)a=0$ readily gives $a^2(bc)=(a^2b)c$. The identities $(a^2b)c=(a^2c)b,(a^2b)b=(b^2a)a$ and $(a^2b)a=a^3b$ are valid for any commutative power-associative algebra over a field of characteristic two and we refer to \cite{5} for a proof.
	\begin{theorem}\label{peircejordan2}
		Let $A$ be a commutative power-associative $K$-algebra of dimension $m\ge 1$ that is not a nilalgebra. If ${\rm char(K)}=2$ there is a nonzero idempotent element
		$e\in A$ that induces the decomposition $A=A_0\oplus A_1$ where:
		\begin{center}
			$A_0=\{a\in A: ea=0\}$ \\
			$A_1=\{a\in A: ea=a\}$
		\end{center}
		are such that $A_0$ and $A_1$ are orthogonal subspaces. Moreover, if $A$ is a Jordan algebra, then $A_1$ and $A_0$ are subalgebras.
	\end{theorem}
	\begin{proof}
		According to \cite{11}, $A$ contains a nonzero idempotent $e$. Consider $x=ex+(ex+x)$. Notice that $e(ex)=(e^2x)e=e^3x=ex$ and $e(ex+x)=ex+ex=0$ and in fact $ex\in A_1,ex+x\in A_0$. Given $a_i\in A_i$, it follows that $a_1a_0=(ea_1)a_0=(ea_0)a_1=0$. Now suppose that $A$ is Jordan. Given $a_i,b_i\in A_i$ for $i\in\{0,1\}$, the conclusion that $A_0$ is a subalgebra follows by computing $e(a_0b_0)=e^2(a_0b_0)=(e^2a_0)b_0=0$ and in a similar fashion $e^2(a_1b_1)=(e^2a_1)b_1=a_1b_1$, therefore $A_0$ and $A_1$ are orthogonal subalgebras of $A$.
	\end{proof}
	
	\begin{example}
		There are commutative power-associative algebras for which the components in the decomposition provided by Theorem \ref{peircejordan2} are not subalgebras. For instance, consider the commutative algebra $A$ over a field of characteristic two with basis $\{e_1,e_2,e_3,e_4\}$ and nonzero products given by:
		$$ e_1^2=e_1,e_1e_2=e_2,e_1e_3=e_3,e_2e_3=e_4.$$
		
		A direct computation shows that $A$ is power-associative. Following the same notation as in Theorem \ref{peircejordan2}, with respect to $e_1$ we have $A=A_0\oplus A_1$ where $A_1=\langle e_1,e_2,e_3\rangle$ but $e_2e_3=e_4\not \in  A_1$ hence $A_1$ is not a subalgebra.
	\end{example}
	When $A$ has a Peirce decomposition $A=A_1(e)\oplus  A_{1/2}(e)\oplus A_0(e)$, we say that the Peirce decomposition of $A$ with respect to $e$ has dimensions $(a_1,a_{1/2},a_0)$ where $a_i=\dim A_i$ for $i\in\{1,\tfrac 12,0\}$ to simplify our notation.
	
	Throughout the text, $A$ may contain more than one idempotent hence other Peirce decompositions may be valid for the same algebra. Moreover, when $A$ is not a commutative associative algebra its Peirce decomposition may have an extra component $A_{1/2}$ with respect to the eigenvalue $1/2$. In our analysis we choose an idempotent $e$ such that dim$(A_1 (e))$ is maximal to give an ordering for the cases, that is, if $A$ has idempotents $e$ and $e'$ such that $A$ has type $(2,1)_e$ and $(3,0)_{e'}$ we only consider the decomposition with respect to $e'$, if more than one idempotent satisfies this condition, we take one that maximizes dim$(A_{1/2}(e))$ when applicable, or else we choose one arbitrarily. This will be our main approach when dealing with types and every time an algebra may be presented by a type which was already studied we shall mention before interrupting the analysis.
	
	\section{Commutative and Associative algebras}
	
	\subsection{Dimension one}
	
	Algebras of dimension 1 may only assume two forms up to isomorphism: $J_{1,1}=\langle n\mid\rangle$
	and $J_{1,2}=\langle f\mid f^2=\alpha f\rangle$ for $\alpha\neq0$ and by the change of basis to $\{\frac{f}{\alpha}\}$ we consider $J_{1,2}\cong\langle e\mid e^2=e\rangle$ and we only take the latter form.
	\subsection{Dimension two} 
	Let $A$ be a commutative associative algebra of dimension two, then $A$ is isomorphic to one of the following algebras:
	\begin{center}
		$J_{2,1}=\lan a,b \mid\ran$; \hspace{0,5cm}
		$J_{2,2}=\lan a,a^2\mid \ran$; \hspace{0,5cm}
		$J_{2,3}=\lan e,a\mid e^2=e \ran$; \hspace{0,5cm}\\ 
		$J_{2,4}^\rho=\lan 1,a\mid a^2=\rho1\ran$, $\rho\in K$;  \hspace{0,5cm} 
		$J_{2,5}^{\kappa}=\lan 1, a\mid a^2=\kappa 1+a\ran,\kappa\in K$.
	\end{center}
	In Section \ref{isomorphismdimension2} we present a necessary and sufficient condition for isomorphism between these families of algebras.\\
	
	For the following, $A$ denotes a two-dimensional algebra.
	
	\subsubsection{Nilalgebras:} Assume that $A$ is a nilalgebra. We may use Corollary \ref{nilbasis3} to conclude that its nilindex is limited by 3 and since the only algebra of nilindex 1 has dimension 0, we check for nilindexes 2 and 3. If $A$ has nilindex 3, regardless of the underlying field, Lemma \ref{nilbasis2} can be applied to produce the algebra described by table $J_{2,2}$. If nil$(A)=2$ a basis $\{x,y\}$ for $A$ satisfies $y^2=x^2=0,xy=\alpha x + \beta y$ for $\alpha,\beta\in K$. Notice that: 
	
	$$
	\begin{cases}
		0=x^2y=\beta(xy)\implies xy=0 \text{ or } \beta=0\\	
		0=xy^2=\alpha (xy) \implies xy=0 \text{ or } \alpha=0
	\end{cases}
	$$
	In any case $xy=0$, therefore this is the zero algebra described by $J_{2,1}$.\\
	
	If $A$ is not a nilalgebra, it contains a nonzero idempotent $e$ by Proposition \ref{idempotent}, hence Lemma \ref{newpeirce} applies. Using the decomposition $A=A_1(e)+A_0(e)$ with ${\rm dim}(A_1(e))\geq 1$ the types to be considered are $(2,0)$ (which is unital) and type $(1,1)$ .
	
	\subsubsection{Type (2,0)}
	In this case, $e$ is the unity element and will be denoted by $1$. Take a basis $\{1,a\}$, we have:
	\begin{equation}\label{equation1}
		a^2=\rho1+\theta a,\ {\rm for\ some} \ \rho,\theta\in K.
	\end{equation}
	
	If $\theta=0$ in equation \ref{equation1} we immediately have algebra $J_{2,4}^\rho$. For $\theta \neq 0$, a change of basis to $\{1,\theta^{-1}a\}$ gives an isomorphism $A(\rho,\theta)\cong J_{2,5}^{\rho\theta^{-2}}$. Reciprocally, any algebra in the form $J_{2,5}^{\kappa}$ is isomorphic to $A(\kappa,1)$. Therefore, we have a correspondence between algebras in the form $A(\rho,\theta)$ and algebras $J_{2,5}^\kappa$.
	
	\subsubsection{ Type (1,1)} 
	Note that $A_0\cong J_{1,1}$ or $A_0\cong J_{1,2}$ since it is one-dimensional. If $A_0\cong J_{1,1}$ we readily have  $A\cong J_{2,3}$. If $A_0\cong J_{1,2}$ it follows that $A\cong\lan e,f\mid e^2=e,f^2=f \ran$, in this case, using a change of basis to $\{e+f,f\}$ 
	we have $A\cong J^0_{2,5}$.
	\subsection{Isomorphism classes}\label{isomorphismdimension2}
	Using the following table, it suffices to check for isomorphisms between pairs of the form $(J_{2,4}^\rho, J_{2,5}^\kappa)$, $(J_{2,5}^\kappa,J_{2,5}^k)$ 
	and $(J_{2,4}^\rho,J_{2,4}^r)$.
	\begin{table}[!h]
		\begin{adjustbox}{max totalsize={\textwidth}{0.88\textheight},center}
			\begin{tabular}{|l|l|}
				\hline
				Algebra &  Notes \\ \hline
				$J_{2,1}$ & Nilpotent of nilindex 2 \\ \hline
				$J_{2,2}$ & Nilpotent of nilindex 3\\ \hline
				$J_{2,3}$ & Nonunital and contains an idempotent\\ \hline
				$J^\rho_{2,4}$ & Unital \\ \hline
				$J^\kappa_{2,5}$ &  Unital \\ \hline
			\end{tabular}
		\end{adjustbox}
	\end{table}
	
	\subsubsection{Isomorphism over fields of characteristic not $2$}
	
	Over fields of $char\neq 2$ algebras $ J_{2,5}^\kappa$ and $J^\rho_{2,4}$ for $\rho=\kappa+\frac{1}{4}$ are isomorphic using a change of basis to $\{1,a-1/2\}$. Therefore it suffices to check for isomorphisms in the case $(J_{2,4}^\rho,J_{2,4}^r)$ for $r\neq \rho$. The method presented below will be extensively used and as such it is necessary to make a complete first example.\\
	
	\begin{proposition} There is an isomorphism $J_{2,4}^\rho\cong J_{2,4}^r$ if and only if there exists $\beta \in K$ with $\beta\neq 0$ and $\rho=\beta^2r$.
	\end{proposition}
	\begin{proof}
		In fact, suppose that there is an isomorphism
		$$h\colon J_{2,4}^\rho=\lan 1,a\ |\ a^2=\rho1 \ran \to J_{2,4}^r=\lan 1,x\ | \ x^2=r1\ran$$
		
		where we used some abuse of notation in writing $1$ for both unities. That is, $h$ has to satisfy:
		
		\begin{enumerate}
			\item $h(1)=1$ and $h(a)=\alpha1+\beta x$ for some $\alpha,\beta \in K$
			\item $\det ([h]) = \left|\begin{array}{cc}
				1 & \alpha  \\
				0 & \beta 
			\end{array}\right|=\beta\neq 0 $
			\item $h(a)^2=h(a^2)=\rho h(1)=\rho 1$;
		\end{enumerate}
		
		Using conditions (1) and (3) we may write $(\alpha^2+\beta^2r)1+2\alpha\beta x=\rho1$ from where we deduce:
		$$
		\begin{cases}
			\beta\neq 0\hspace{0,5cm} &(By\ Condition\ (2))\\
			2\alpha\beta=0 &\implies \alpha=0 \\
			\rho=\alpha^2+\beta^2r &\implies \rho=\beta^2r 
		\end{cases}
		$$
		
		This proves the implication. For the converse,
		suppose such a $\beta$ exists. Consider the basis $\{1,\beta x\}$ for $J_{2,4}^r$. Then $(\beta x)^2=(\beta^2r)1=\rho1$, this shows that $J_{2,4}^r\cong J_{2,4}^{\beta^2r}=J_{2,4}^{\rho}$.
	\end{proof}
	
	\subsubsection{Isomorphism over fields of characteristic two.} 
	\leavevmode\par
	
	\begin{proposition} Algebras of the family $J_{2,4}^\rho$ and $J_{2,5}^{\kappa}$ are not isomorphic.
	\end{proposition}
	\begin{proof}
		Suppose there is an isomorphism 
		$$h\colon J_{2,4}^\rho=\lan 1,a\mid a^2=\rho 1 \ran \to J_{2,5}^{\kappa}=\lan 1,x \mid x^2=\kappa 1+x\ran$$
		where we used the same abuse of notation in writing $1$ for both units, this means that $h$ has to satisfy:
		\begin{enumerate}
			\item $h(1)=1$ and $h(a)=\alpha1+\beta x$ for some $\alpha,\beta \in K$
			\item $\det ([h]) = \left|\begin{array}{cc}
				1 & \alpha  \\
				0 & \beta 
			\end{array}\right|=\beta\neq 0 $
			\item $h(a)^2=h(a^2)=\rho h(1)=\rho 1$;
		\end{enumerate}
		then $(\alpha^2+\beta^2\kappa)1+(2\alpha\beta+\beta^2) x=(\alpha^2+\beta^2\kappa)1+\beta^2 x= \rho 1 $ which implies $\beta=0$, a contradiction.
	\end{proof}
	
	
	\begin{proposition} There is an isomorphism $J_{2,4}^\rho\cong J_{2,4}^r$ if and only if there exists $(\alpha,\beta)\in K^2$ with $\beta \neq 0$ satisfying $\alpha^2+\beta^2r=\rho.$
	\end{proposition}
	\begin{proof}
		Suppose there is an isomorphism 
		$$h\colon J_{2,4}^\rho=\lan 1,a \mid a^2=\rho1 \ran \to J_{2,4}^r=\lan 1,x \mid x^2=r1\ran$$
		where we used the same abuse of notation in writing $1$ for both unities. The function $h$ has to satisfy:
		\begin{enumerate}
			\item $h(1)=1$ and $h(a)=\alpha1+\beta x$ for some $\alpha,\beta\in K$;
			\item $\det ([h]) = \left|\begin{array}{cc}
				1 & \alpha  \\
				0 & \beta 
			\end{array}\right|=\beta\neq 0 $
			\item $h(a)^2=h(a^2)=\rho h(1)=\rho 1$; 
		\end{enumerate}
		then $(\alpha^2+\beta^2r)1=\rho1$.
		This proves the implication.
		Suppose a pair $(\alpha,\beta)$ satisfying the above conditions exists. Consider the basis $\{1,\alpha1+\beta x\}$ for $J_{2,4}^r$, note that $(\alpha1+\beta x)^2=(\alpha^2+\beta^2r)1=\rho1$. This shows that $J_{2,4}^r\cong J_{2,4}^{\alpha^2+\beta^2r}=J_{2,4}^{\rho}$.\end{proof}
	
	\begin{proposition}  There is an isomorphism $J_{2,5}^{k}\cong J_{2,5}^{\kappa}$ if and only if there exists $\alpha\in K$ such that $\alpha^2+k=\kappa+\alpha$.
	\end{proposition}
	\begin{proof}
		Suppose there is an isomorphism $$h\colon J_{2,5}^{\kappa}=\lan 1,a\mid a^2=\kappa1+a \ran \to J_{2,5}^k=\lan 1,x \mid x^2=k1+x\ran$$
		where we used the same abuse of notation in writing $1$ for both units, this means that $h$ has to satisfy:
		\begin{enumerate}
			\item $h(1)=1$ and $h(a)=\alpha1+\beta x$ for some $\alpha,\beta\in K$;
			\item $\det ([h]) = \left|\begin{array}{cc}
				1 & \alpha  \\
				0 & \beta 
			\end{array}\right|=\beta\neq 0 $	
			\item $h(a)^2=h(a^2)=\kappa h(1)+h(a)$;
		\end{enumerate}
		then  $(\alpha^2+\beta^2k)1+\beta^2x=(\kappa+\alpha)1+\beta x$ from where we deduce:
		\begin{center}
			$
			\begin{cases}
				\beta\neq 0 &\hspace{0,5cm} (By\ Condition\ $(2)$)\\
				\beta^2=\beta	&\implies \beta=1\\
				\alpha^2+\beta^2k=\kappa+\alpha &\implies \alpha^2+k=\kappa+\alpha
			\end{cases}
			$
		\end{center}
		This proves the implication. Conversely, suppose that there exists $\alpha\in K$ such that $\alpha^2+k=\kappa+\alpha$. Note that $(\alpha1+x)^2=\alpha^21+k1+x=(\alpha^2+k)1+x=(\kappa+\alpha)1+x=\kappa1+(\alpha1+x)$. That is, a change of basis to $\{1,\alpha1+x\}$ proves that $J_{2,5}^k\cong\lan 1, \alpha1+x \mid(\alpha1+x)^2=\kappa1+(\alpha1+x)\ran=J_{2,5}^\kappa$.
	\end{proof}
	\subsection{Dimension three.}
	Every commutative associative nilalgebra is nilpotent \cite{11} and since associative nilpotent algebras of dimension at most 4 have been fully classified \cite{9} this case is settled. We listed those algebras below:
	\begin{center}
		$J_{3,1} = \langle a, b, c \mid  \rangle$; \hspace{0.5cm}
		$J_{3,2} = \langle a, b, c \mid a^2 = b \rangle$; \\
		$J_{3,3}^{\delta} = \langle a, b, c \mid a^2 = c, b^2 = \delta c \rangle$, $\delta \neq 0$; \\
		$J_{3,4} = \langle a, b, c \mid a^2 = b, ab = c \rangle$. 
	\end{center}
	The table $A=\lan a,b,c\mid ab=-c,ba=c\ran$ is noncommutative when ${\rm char}(K)\neq2$. In characteristic two it is commutative and associative; we defer this table to the zeropotent algebras in Section \ref{dimension3jordan}, where it occurs as $ZA_3$ \cite{19}. For easier indexation we list the algebras obtained in the following cases below.
	\begin{center}
		$J_{3,5}=\lan 1,a,b\mid  b^2=b \ran $; \hspace{0,5cm}
		$J_{3,6}=\lan 1,a,b \ran $; \hspace{0,5cm}
		$J_{3,7}=\lan 1,a,a^2\ran $\\
		$J_{3,8}=\langle a,b,c\mid a^2=b,c^2=c\rangle$; 
		\hspace{0,5cm}
		$J_{3,9}=\langle a,b,c\mid a^2=a\rangle$; \\
		$J^\rho_{3,10}=\langle a,b,c\mid a^2=a,ab=b,b^2=\rho a\rangle$, $\rho\in K$;\\
		$J^\kappa_{3,11}=\langle a,b,c\mid a^2=a,ab=b,b^2=b+\kappa a\rangle$,\ $\kappa\in K$;\\
		$J^\rho_{3,12}=\langle a,b,c\mid a^2=a,b^2=b,bc=c,c^2=\rho b\rangle$, $\rho \in K$;\\
		$J^\kappa_{3,13}=\langle a,b,c\mid a^2=a,b^2=b+\kappa a,c^2=c,ab=b\rangle$, $\kappa \in K$;\\
		$J_{3,14}=K(\alpha)$, an extension of $K$ of degree 3.
	\end{center}
	
	Suppose that $A$ is not a nilalgebra. By Proposition \ref{idempotent}, $A$ contains an idempotent $e\neq 0$ and Lemma \ref{newpeirce} applies.  Since ${\rm dim }A_1(e)\geq 1$, it suffices to consider types $(3,0),(2,1)$ and $(1,2)$.
	
	\subsubsection{Type $(3,0)$.}
	We study the cases ${\rm dim}({\rm Rad}(A))\in \{0,1,2\}$.\\
	
	\textit{Case ${\rm dim}({\rm Rad}(A))=0$:} Since $A/{\rm Rad}(A)$ is semisimple and in this case $A\cong A/{\rm Rad}(A)$, we may apply the Wedderburn-Artin Theorem to show that there exist positive integers $n_1, \ldots, n_r$ and division algebras $D_1, \ldots, D_r$ such that $A/J(A) \cong M_{n_1}(D_1) \oplus \cdots \oplus M_{n_r}(D_r)$, where $M_{n_i}(D_i)$ is a matrix algebra with coefficients in $D_i$, for every $1 \leq i \leq r$.
	Observe that since $0={\rm Rad}(A)=J(A)$, and $A$ is a $K$-algebra, each $D_i$ is also an algebra over $K$. Furthermore, since any algebra of matrices of order greater than 1 has dimension at least four over $K$, we conclude that $A$ is either:
	
	\begin{itemize}
		\item A direct sum of $K\oplus K\oplus K$; In this case, the sum of the identity elements of each factor is the identity element of $A$. Taking a basis of orthogonal idempotents $\{e_1,e_2,e_3\}$ and applying a change of basis to  
		$\{e_1 + e_2, e_2, e_3\}$, we obtain $K\oplus K\oplus K \cong J^0_{3,13}$.
		\item A direct sum of $K$ with a division algebra over $K$ of dimension 2, here denoted by $K'$. By the description of commutative associative algebras of dimension 2, either $K' \cong J^\rho_{2,4}$ for $\rho$ not a square or $K' \cong J^\kappa_{2,5}$  for $\kappa$ such that $p(X)=X^2-X-\kappa$ is irreducible over $K$; under these conditions $K\oplus K'\cong J^\rho_{3,12}$ or $K\oplus K'\cong J^\kappa_{3,13}$.
		Since $K'$ is associative, to be a division algebra is to be free of zero divisors, below we explicitly determine the conditions for the algebras $J_{2,4}^\rho$ and  $J_{2,5}^\kappa$ to satisfy this condition.\\ 
		
		Conditions for $J_{2,4}^\rho$ to be free of zero divisors (over any field):
		Let $x_11+x_2a, x_31+x_4a \in        
		J_{2,4}^\rho \setminus \{ 0 \}$ with $x_1,x_2,x_3,x_4\in K$ be such that $$(x_11+x_2a)(x_31+x_4a)=0.$$ Since $A$ is  associative and unital, no invertible element $x_1+x_2a\in A$ is a zero divisor. Notice that:
		if $(x_1+x_2a)^{-1}$ exists, then 
		$$(x_1+x_2a)^{-1}=\frac{x_1}{x_1^2-\rho x_2^2}1-\frac{x_2}{x_1^2-\rho x_2^2}a$$
		For a fixed nonzero element $x_11+x_2a$, the inverse exists if and only if $x_1^2-\rho x_2^2\neq0$. If $x_2=0$, this holds because $x_1\neq0$. If $x_2\neq0$, it is equivalent to $\rho\neq(x_1/x_2)^2$. Consequently, every nonzero element is invertible if and only if $\rho$ is not a square in $K$. If $\rho=r^2$, the nonzero elements $a-r1$ and $a+r1$ have product zero (and may coincide in characteristic two).
		
		Conditions for $J_{2,5}^\kappa$ to be free of zero divisors over a field of characteristic two: For this case we have: If $(x_1+x_2a)^{-1}$ exists, then
		$$(x_1+x_2a)^{-1}=\frac{-x_1-x_2}{-x_1^2- x_1x_2+x_2^2 \kappa}1+\frac{x_2}{-x_1^2- x_1x_2+x_2^2 \kappa}a$$
		
		For a fixed nonzero element $x_11+x_2a$, the inverse exists if and only if $x_1^2+x_1x_2-\kappa x_2^2\neq0$. The case $x_2=0$ is automatic. If $x_2\neq0$, in characteristic two this is equivalent to $\kappa\neq(x_1/x_2)^2+x_1/x_2$. Therefore $J_{2,5}^{\kappa}$ is a division algebra if and only if $\kappa\notin\{r^2+r:r\in K\}$, equivalently, $X^2-X-\kappa$ is irreducible over $K$.
		
		\item An extension of degree 3 of $K$. This case was denoted by $J_{3,14}$.
	\end{itemize}
	
	Case ${\rm dim}({\rm Rad}(A))=1$: Let $0\neq a\in {\rm Rad}(A)$. 
	Since ${\rm Rad}(A)$ is a nilideal of $A$, we can prove $a^2=0$. Take a basis $\{1,a,b\}$ of $A$ satisfying $a^2=0,ab=p_1a,b^2=\gamma_11+\gamma_2a+\gamma_3b$. Notice that, if $p_1 \neq 0$, then a  change of basis 
	$b\mapsto b-p_1 1$ 
	provides a new basis that satisfies $a^2=ab=0$. Abusing notation if needed, the basis is taken to satisfy this extra condition. The equality $0=(ab)b=ab^2=\gamma_1a$ implies $\gamma_1=0$. Furthermore, we state that $\gamma_3 \neq 0$. In fact, $\gamma_3 = 0$ implies $b \in {\rm Rad}(A)$  which is a contradiction, since $\{a,b\}$ is linearly independent and  ${\rm dim}({\rm Rad}(A))=1$. Finally, a change of basis $b \mapsto \gamma_3^{-1}(\gamma_3^{-1}\gamma_2a+b)$ results in the algebra $J_{3,5}$.\\
	
	Case ${\rm dim}({\rm Rad}(A))=2$: For ${\rm nil}({\rm Rad}(A)) = 2$, let $\{a, b\}$ be a basis of ${\rm Rad}(A)$ satisfying $a^2 = b^2 = 0$ and $ab=\rho_1a+\rho_2b$, $\rho_1, \rho_2 \in K$. Since $0=ab^2=(ab)b$ and $0=a^2b=a(ab)$, the equalities $\rho_1=\rho_2=0$ hold. Therefore, this is a zero algebra adjoined to 1, as described by 
	$J_{3,6}$.  
	For ${\rm nil}({\rm Rad}(A)) = 3$, by Lemma \ref{nilbasis2}, a basis $\{a,a^2\}$ for ${\rm Rad}(A)$ may be taken to satisfy $a^3=0$. This automatically gives the algebra $J_{3,7}$.
	
	\subsubsection{Type $(2,1)$.} 
	Since $A_1(e)$ is a two-dimensional unital algebra with identity element $e$, we have $A_1(e)\cong J_{2,4}^\rho$ or $A_1(e)\cong J_{2,5}^\kappa$. If $f\in A_0(e)$ is a nonzero idempotent, then $e+f$ is the identity element of $A$, violating the maximality of $\dim A_1(e)<3$. Hence $A_0(e)\cong J_{1,1}$. Therefore, the resulting algebras are isomorphic to $J_{2,4}^\rho\oplus J_{1,1}$ or $J_{2,5}^\kappa\oplus J_{1,1}$ which are isomorphic to, respectively, $J^\rho_{3,10}$ and $J^\kappa_{3,11}$.
	
	\subsubsection{Type $(1,2)$.} In this case $\lan e\ran=A_1(e)$ is a unital one-dimensional algebra, we conclude that $A_1(e)\cong J_{1,2}$. If $f\in A_0(e)$ is a nonzero idempotent, then $e+f$ is a nonzero idempotent satisfying $\lan e,f\ran\subseteq A_1(e+f)$, violating the maximality of $\dim A_1(e)=1$. Thus, $A_0(e)$ is isomorphic to either $J_{2,2}$ or $J_{2,1}$. For $A_0(e)\cong J_{2,2}$ we have $A\cong J_{1,2}\oplus J_{2,2}\cong J_{3,8}$ and for  $A_0(e)\cong J_{2,1}$ we have $A\cong J_{1,2}\oplus J_{2,1}\cong J_{3,9}$.
	
	\subsubsection{Isomorphism classes}
	\begin{table}[!h]
		\centering
		\begin{tabular}{|l|l|}
			\hline
			Algebra & Note \\ \hline
			$J_{3,1}$ & Zero algebra \\ \hline
			$J_{3,2}$ & Nilpotent \\ \hline
			$J_{3,3}^{\delta}$ & Nilpotent \\ \hline
			$J_{3,4}$ & Nilpotent \\ \hline
			$J_{3,5}$ 
			& Unital, $\dim \operatorname{Rad}(J_{3,5})=1$ \\ \hline
			$J_{3,6}$ 
			& Unital, $\dim \operatorname{Rad}(J_{3,6})=2$, 
			$\operatorname{nil}(\operatorname{Rad}(J_{3,6}))=2$ \\ \hline
			$J_{3,7}$ 
			& Unital, $\dim \operatorname{Rad}(J_{3,7})=2$, 
			$\operatorname{nil}(\operatorname{Rad}(J_{3,7}))=3$ \\ \hline
			$J_{3,8}$ 
			& Nonunital, $J_{2,2}\oplus J_{1,2}$, 
			$\operatorname{nil}(\operatorname{Rad}(J_{3,8}))=3$ \\ \hline
			$J_{3,9}$ 
			& Nonunital, $J_{2,1}\oplus J_{1,2}$, 
			$\operatorname{nil}(\operatorname{Rad}(J_{3,9}))=2$ \\ \hline
			$J_{3,10}^{\rho}$ 
			& Nonunital, $J_{2,4}^{\rho}\oplus J_{1,1}$, 
			$\operatorname{nil}(\operatorname{Rad}(J_{3,10}^{\rho}))=2$ \\ \hline
			$J_{3,11}^{\kappa}$ 
			& Nonunital, $J_{2,5}^{\kappa}\oplus J_{1,1}$, 
			$\operatorname{nil}(\operatorname{Rad}(J_{3,11}^{\kappa}))=2$ \\ \hline
			$J_{3,12}^{\rho}$ 
			& Unital, $J_{2,4}^{\rho}\oplus J_{1,2}$ \\ \hline
			$J_{3,13}^{\kappa}$ 
			& Unital, $J_{2,5}^{\kappa}\oplus J_{1,2}$ \\ \hline
			$J_{3,14}$ 
			& Field \\ \hline
		\end{tabular}
	\end{table}
	In virtue of the table above, it remains to check the following isomorphism cases.\\
	
	Since $J_{3,5}\cong J_{2,4}^0\oplus J_{1,2}$ and $J_{3,12}^{\rho}\cong J_{2,4}^{\rho}\oplus J_{1,2}$, their isomorphism criterion is
	$$
	J_{3,5}\cong J_{3,12}^{\rho}\quad\Longleftrightarrow\quad
	\begin{cases}
		\rho=0,&{\rm char}(K)\neq2,\\
		\rho\text{ is a square in }K,&{\rm char}(K)=2.
	\end{cases}
	$$
	Indeed, the quadratic factor has a nonzero nilpotent element precisely when its defining polynomial has a repeated root in $K$. In characteristic two, if $\rho=r^2$, replacing $c$ by $c+rb$ in $J_{3,12}^{\rho}$ gives $J_{3,12}^0$, and the basis $\{a+b,c+rb,a\}$ gives $J_{3,5}$. In characteristic different from two a repeated root of $X^2-\rho$ occurs precisely when $\rho=0$. Similarly,
	$$
	J_{3,5}\cong J_{3,13}^{\kappa}\quad\Longleftrightarrow\quad
	{\rm char}(K)\neq2\ \text{ and }\ \kappa=-\tfrac14.
	$$
	The polynomial $X^2-X-\kappa$ has a repeated root precisely in this case.\\
	
	We conclude $J_{3,6}\not\cong J_{3,12}^\rho$ and $J_{3,6}\not\cong J_{3,13}^\kappa$ by counting their nonzero idempotents, for the former algebra, there is only one while for the latter at least two.\\
	
	Since $J_{3,7}$ does not have an idempotent element that induces a Peirce decomposition satisfying $\dim A_0>0$ we conclude that it is not isomorphic to $J^\rho_{3,12},J^\kappa_{3,13}$\\
	
	Notice that $\dim Ann(J_{3,9})=2, \dim Ann(J_{3,10}^\rho)=1,\dim Ann(J_{3,11}^\kappa)=1$ hence there is no isomorphism between $J_{3,9}$ and any of these two algebras.\\
	
	If ${\rm char}(K)\neq 2$, then
	$J_{2,4}^{\rho}\cong J_{2,5}^{\kappa}$ if and only if there exists
	$s\in K\setminus\{0\}$ such that
	$$
	\rho=s^2\left(\kappa+\frac14\right).
	$$
	Using the direct-sum descriptions of the three-dimensional algebras,
	the same necessary and sufficient condition characterizes the
	isomorphisms
	$$
	J_{3,10}^{\rho}\cong J_{3,11}^{\kappa}
	\qquad\text{and}\qquad
	J_{3,12}^{\rho}\cong J_{3,13}^{\kappa}.
	$$
	If ${\rm char}(K)=2$, the two algebras in each of these three pairs
	are nonisomorphic for every choice of $\rho,\kappa\in K$.

	Within the family $J_{3,3}^{\delta}$, where $\delta\neq0$, we have
	$$
	J_{3,3}^{\delta}\cong J_{3,3}^{\epsilon}
	\quad\Longleftrightarrow\quad
	\delta/\epsilon\text{ is a square in }K.
	$$
	This was obtained by Graaf \cite{9}. The within-family criteria for
	$J_{3,10}^{\rho},J_{3,11}^{\kappa},J_{3,12}^{\rho},J_{3,13}^{\kappa}$
	are the corresponding two-dimensional criteria: one recovers the
	quadratic factor by taking the quotient by the annihilator in the
	nonunital cases, and by uniqueness of the decomposition into local
	factors in the unital cases. The label $J_{3,14}$ denotes the cubic
	field extensions up to $K$-algebra isomorphism.
	
	\section{Jordan algebras}
	
	Using central extensions and the Skjelbred-Sund method to construct nilpotent algebras over arbitrary fields, Hegazi and Abdelwahab \cite{14} classified nilpotent Jordan algebras of dimension at most 3. Incidentally, every nilpotent Jordan algebra of dimension at most 3 is associative.
	
	\subsection{Dimension one.}
	Every one-dimensional commutative algebra is associative, so this variety was precisely described in the last sections.
	
	\subsection{Dimension two.}
	\subsubsection{The case ${\rm char}(K)\neq 2$.}
	Two-dimensional Jordan algebras have been described by Diniz et al. \cite{8}. 
	The associative representatives are 
	$J_{2,1},J_{2,2}$, $J_{2,3}$ and $J^\rho_{2,4}$ and the unique non-associative algebra is
	\begin{center}
		$J_{2,6}=\lan e,a\ |\ e^2=e,ea=\frac{a}{2}\ran$.
	\end{center}
	
	\subsubsection{The case ${\rm char}(K)=2$.}
	\begin{center}
		$J_{2,7}=\lan a,b\ |\ ab=a\ran$. 
	\end{center}
	
	Nilalgebras. For nilalgebras of nilindex $3$, Lemma \ref{nilbasis2} applies from where algebra $J_{2,2}$ is derived. Suppose ${\rm nil}(A)=2$ and fix a basis $\{a,b\}$ with $ab=\rho_1a+\rho_2b$. If $ab=0$ we have the zero algebra $J_{2,1}$, so we suppose $ab\neq 0$. If $\rho_1=0$, a change of basis to $\{a',b'\}$ where $a'=b$ and $b'=\rho_2^{-1}a$ gives $(a')^2=(b')^2=0$ and $a'b'=a'$, this is algebra $J_{2,7}$. If $\rho_1\neq 0$, we consider a change of basis to $\{a',b'\}$, where $a'=\rho_1a+\rho_2b$ and  $b'=\rho_1^{-1}b$, and again $(a')^2=2\rho_1\rho_2a'=0$, $(b')^2=0$ and $a'b'=a'$. In both cases we have algebra $J_{2,7}= \lan a,b\ |\ ab=a\ran $. Note that $J_{2,7}$ is a non-nilpotent Jordan nilalgebra, since $(\cdots ((ab)b) \cdots )b=a\neq 0$.  
	
	\subsubsection{Algebras that are not nilalgebras} In this case, all algebras are associative. In fact,
	since these algebras are power-associative we may consider their Peirce decomposition as in the previous sections.
	For type $(2,0)$, there exists $a \in A$ such that $\{1,a\}$ is a basis of $A$. Write $a^2=\rho_11+\rho_2a$, for some $\rho_1, \rho_2 \in K$. If $\rho_1\neq 0$ the set $\{a,a ^2\}$ is linearly independent, therefore $A=\lan a,a^2 \ran$, hence $A$ is  associative. Otherwise, suppose $a^2=\rho_2a$. If $\rho_2=0$, then $A =\langle a+1,(a+1)^2 \rangle$ and the same conclusion follows. If $\rho_2 \neq 0$, by a change of basis to $a\mapsto \rho_2^{-1}a$ we may further consider $\rho_2=1$. 
	In this case, the set $\{t+a,(t+a)^2\}$ is linearly independent for any $t\notin \{0,1\}$ and $A$ is again associative. Hence, every Jordan algebra of type (2,0) is associative. Algebras of type $(1,1)$ are a direct sum of orthogonal commutative associative algebras. Thereby, all algebras obtained have already been studied in the previous section.
	\subsection{Dimension three}\label{dimension3jordan}
	\begin{center}
		$J_{3,15}^{\ell}=\lan a,a^2,b\mid ab=b,\ b^2=\ell a^2\ran,\quad \ell\in K$.
	\end{center}
	
	\subsubsection{Nilalgebras.}
	Over fields of characteristic different from two, every Jordan nilalgebra is nilpotent and the results in \cite{14} prove that $A$ is associative. Thus, it suffices to study the cases where the field has characteristic two.\\
	
	Any three-dimensional power-associative nilalgebra of nilindex $4$ is associative, by Lemma \ref{nilbasis2}. Suppose now that $A$ is a three-dimensional Jordan nilalgebra of nilindex
	$3$ over a field $K$ of characteristic $2$ with more than two elements.
	Thus $z^3=0$ for every $z\in A$, and there exists an element $x\in A$
	such that $x^2\neq0$. The vectors $x$ and $x^2$ are linearly independent.
	Indeed, if $x^2=\delta x$ for some $\delta\in K$, then
	$0=x^3=\delta x^2$, so $\delta=0$, contradicting $x^2\neq0$.
	We therefore complete them to a basis $\{x,x^2,y\}$ of $A$ and write
	$$
	xy=\rho_1x+\rho_2x^2+\rho_3y,
	\qquad
	y^2=\alpha_1x+\alpha_2x^2+\alpha_3y,
	\qquad \rho_i,\alpha_i\in K.
	$$
	
	We first show that every square annihilates $A$. Given $z,w\in A$ and
	$t\in K$, commutativity and characteristic $2$ give
	$$
	0=(z+tw)^3
	=(z^2+t^2w^2)(z+tw)
	=t z^2w+t^2w^2z.
	$$
	Taking $t=1$ yields $z^2w+w^2z=0$, hence $z^2w=w^2z$.
	Choosing $t\notin\{0,1\}$ now gives $t(1+t)z^2w=0$, so $z^2w=0$.
	In particular, $x^2A=0$ and $y^2A=0$.\\
	
	We next prove that $y^2\in\langle x^2\rangle$. Expanding $y^2x=0$
	in the chosen basis, we obtain
	$$
	0=\alpha_3\rho_1x
	+(\alpha_1+\alpha_3\rho_2)x^2
	+\alpha_3\rho_3y.
	$$
	Therefore
	$$
	\alpha_3\rho_1=0,
	\qquad \alpha_1+\rho_2\alpha_3=0,
	\qquad \alpha_3\rho_3=0.
	$$
	Suppose that $\alpha_3\neq0$. The first and third equations listed above then imply
	$\rho_1=\rho_3=0$, so $xy=\rho_2x^2$. Consequently,
	\begin{align*}
		0=y^3
		&=\alpha_1xy+\alpha_2x^2y+\alpha_3y^2\\
		&=\alpha_1\rho_2x^2
		+\alpha_3(\alpha_1x+\alpha_2x^2+\alpha_3y)\\
		&=\alpha_1\alpha_3x
		+(\alpha_1\rho_2+\alpha_2\alpha_3)x^2
		+\alpha_3^2y.
	\end{align*}
	The coefficient of $y$ must vanish, giving $\alpha_3^2=0$, a
	contradiction. Hence $\alpha_3=0$. The equation
	$\alpha_1+\rho_2\alpha_3=0$ now gives $\alpha_1=0$, and we conclude that $y^2=\alpha_2x^2$.\\
	
	Recall that $xy=\rho_1x+\rho_2x^2+\rho_3y$ with $\rho_1,\rho_2,\rho_3\in K$. Assuming $\rho_3\neq 0$ the change of basis $\{{\rho_3}^{-1}x,{\rho_3}^{-2}x^2,\rho_1{\rho_3}^{-1}x+(\rho_1{\rho_3}^{-2}+\rho_2{\rho_3}^{-1})x^2+y\}$ gives the algebra $J^\ell_{3,15}$ with $\ell=\rho_1^2+\alpha_2\rho_3^2$. If $\rho_3= 0$, $\rho_1\neq 0$ and $\alpha_2 \neq 0$ the change of basis $\{{\rho_1}^{-1}y,\alpha_2 {\rho_1}^{-2}x^2,x+\rho_2{\rho_1}^{-1}x^2\}$ gives the algebra $J^\ell_{3,15}$ with $\ell=\rho_1^2/\alpha_2\neq0$. If $\rho_3= 0$, $\rho_1\neq 0$ and $\alpha_2 = 0$ the change of basis $\{x+\rho_2{\rho_1}^{-1}x^2+{\rho_1}^{-1}y,x^2,x+(1+\rho_2{\rho_1}^{-1})x^2\}$ gives the algebra $J_{3,15}^1$. For $\rho_3=\rho_1=0$ we have an associative algebra, since $A^2 = \langle x^2\rangle$ and $x^2 A = 0$.\\

	For nilalgebras of nilindex two over fields of characteristic two, the classification was completely done by the work of Cedilnik and Jerman \cite{19} and we list them below.\\
	$ZA_1=\lan a,b,c\ran;\hspace{0,3cm}
	ZA_2(0)=\lan a,b,c\ |\ ac=a,bc=a\ran;\hspace{0,3cm}
	ZA_3=\lan a,b,c\ |\ ba=c\ran;$\\
	$ZA_4=\lan a,b,c\ |\ ac=b,bc=c\ran;\hspace{0,3cm}
	ZA_5=\lan a,b,c\ |\ ac=a,bc=c\ran;$\\
	$ZA_6=\lan a,b,c\ |\ ac=a,bc=b\ran; \hspace{0,3cm}
	ZA_2(\lambda\neq 0)=\lan a,b,c\ |\ ac=a+\lambda b,bc=a\ran;$\\
	$ZA_7=\lan a,b,c\ |\ ac=b,bc=a\ran;\hspace{0,3cm} ZA_8(\lambda)=\lan a,b,c \ | \ ac=\lambda b,bc=a\ran;\\
	ZA_9=\lan a,b,c \ | \ ab=a,ac=b,bc=c\ran;$\\
	$ZA_{10}=\lan a,b,c \ | \ ab=a,ac=b,bc=\lambda a+c\ran,\lambda \ not \ a \ square;$\\
	$ZA_{11}(\lambda)=\lan a,b,c \ | \ ab=a,ac=b,bc=\lambda a+b+c\ran, for\ all\ x\in K, \lambda+x+x^2\neq 0;$\\
	$ZA_{14}(\lambda)=\lan a,b,c \ | \ ab=a,ac=c,bc=b+\lambda c\ran;$\\
	$ZA_{15}(\lambda,\mu)=\lan a,b,c \ | \ ab=c,ac=\lambda b,bc=\mu a+b\ran, where\ (0,0,0)\ is\ the\ only\ root\ of\ p(x,y,z)=xy+x^2\mu+y^2\lambda+z^2;$\\
	$ZA_{16}(\lambda,\mu)=\lan a,b,c \mid ab=c,ac=\lambda b, bc=\mu a\ran, where\ (0,0,0)\neq (x,y,z)\ implies\ x^2\mu+y^2\lambda+z^2\neq0$.
	\subsubsection{Isomorphism classes over a field of characteristic two.}
	\leavevmode\par
	
	\begin{proposition} Algebras $J^\lambda_{3,15},J^l_{3,15}$ are isomorphic over a field $K$ of ${\rm char}(K)=2$ if and only if exists $\alpha,\beta \in K,\beta\neq 0$ such that $\alpha^2l+1 \neq 0$ and $l\beta^2= \lambda (1+\alpha^2 l)$.
	\end{proposition}
	\begin{proof}
		Suppose there is an isomorphism
		$$h\colon J^\lambda_{3,15}=\lan a,a^2,b \mid ab=b,b^2=\lambda a^2\ran\to J^l_{3,15}=\lan x,x^2,y \mid xy=y,y^2=l x^2\ran$$
		
		In this case, $h$ must satisfy:
		\begin{enumerate}
			\item $h(a)=\alpha_1x+\alpha_2x^2+\alpha_3y,h(a^2)=(\alpha_1^2+\alpha_3^2l)x^2,h(b)=\beta_1x+\beta_2x^2+\beta_3y$;
			\item  $\det ([h]) = \left|\begin{array}{ccc}
				\alpha_1 & 0 & \beta_1\\
				\alpha_2 & \alpha_1^2+\alpha_3^2l & \beta_2\\
				\alpha_3 & 0 & \beta_3
			\end{array}\right|=(\alpha_1^2+\alpha_3^2l)(\alpha_1\beta_3+\alpha_3\beta_1)\neq 0 $;
			\item $h(ab)=h(a)h(b)=h(b)$ and $h(b^2)=h(b)^2=\lambda h(a)^2$.
		\end{enumerate}
		
		From this we conclude that:
		$$
		\begin{cases}
			\beta_1=0;\\
			\beta_3\neq 0;\\
			\beta_3^2l=\lambda (\alpha_1^2+\alpha_3^2l);\\
			\alpha_1\beta_3=\beta_3;\\
			\alpha_3\beta_3l=\beta_2.
		\end{cases}
		$$
		By Condition (2), $\alpha_3^2l+1\neq 0$. Since $\alpha_1=1$ we have $\beta_3^2l=\lambda (1+\alpha_3^2l)$, proving the implication. For the converse, consider a change of basis to $\{x+\alpha y,(1+\alpha^2l)x^2, \alpha\beta l x^2+\beta y\}$.	
	\end{proof}
	\section{Peirce decomposition and characterization.}
	\subsection{Unital algebras.}   
	For completeness, we make a few introductory remarks. The length of a finite subset $S$ of a $K$-algebra $A$ is defined as follows.
	
	\begin{definition}
		A product of a finite number of elements of $S$ is a word in the letters from $S$. The length of a word equals the number of factors in the corresponding product. If $A$ is unital we consider $1$ as a word of length 0.
	\end{definition}
	
	The set of words with lengths not exceeding $n$ is denoted by $S^n$. For instance, $S^0=\{1\}$ for unital algebras and $S^0=\emptyset$ otherwise. The set $\mathcal{L}_n(S)$ denotes the $K$-vector space generated by $S^n$, additionally $\calL(S)$ stands for the union $\bigcup_{i=0}^{\infty}\mathcal{L}_i(S)$.

	\begin{definition}
		A set $S$ is a generating set of $A$ if and only if $A=\mathcal{L}(S)$ and the length of a generating set $S$ of a finite-dimensional algebra $A$ is $\ell(S)=\min \{k\in\mathbb{Z}_+: \calL_k(S)=A \}$. We also define the length of the $K$-algebra $A$ by $\ell(A)=\max \{\ell(S): \mathcal{L}(S)=A\}$.
	\end{definition}
	We proceed with our analysis of unital algebras over a field of characteristic not $2$. First, notice that every 3-dimensional unital algebra has length limited by 2.
	If $\ell(A)=2$, choose a generating set $S$ with $\mathcal L_1(S)\neq A$. Since $1\in\mathcal L_1(S)$ and $S$ generates the three-dimensional algebra, $\dim\mathcal L_1(S)=2$. Thus $\mathcal L_1(S)=\langle1,x\rangle$ for some $x$, and every element of $S$ belongs to this vector space. Products of two elements of $\langle1,x\rangle$ lie in $\langle1,x,x^2\rangle$. Since $\mathcal L_2(S)=A$, the set $\{1,x,x^2\}$ is a basis. Power-associativity makes the subalgebra generated by $x$ associative; adjoining the identity element to this subalgebra preserves associativity, so $A$ is associative. If $\ell(A)=1$, it follows from \cite[Corollary 3.5]{3} that a basis may be taken to satisfy
	$$
	a^2=\rho1,b^2=\delta1,ab=\gamma1, \qquad \text{for some}\ \rho,\delta,\gamma\in K.
	$$
	Thus, we obtain the following algebras
	\begin{center}
		$ J^{\alpha,\beta}_{3,16}=\lan 1,a,b \mid a^2=\alpha1,\ b^2=\beta1\ran$, $\alpha,\beta\in K$; \\
		$J_{3,17}=\lan 1,a,b\mid ab=1\ran$.
	\end{center}
	$(\rho = 0,\ \delta = 0,\ \gamma = 0):$  Isomorphic to $J^{0,0}_{3,16}$, no change of basis needed;\\
	$(\rho = 0,\ \delta = 0,\ \gamma \neq 0):$ Isomorphic to $J_{3,17}$ changing basis to $\{1,\frac{a}{\gamma},b\}$;\\
	$(\rho = 0,\ \delta \neq 0,\ \gamma = 0):$ Isomorphic to $J_{3,16}^{0,\delta}$, no change of basis performed;\\
	$(\rho = 0,\ \delta \neq 0,\ \gamma \neq 0):$ Isomorphic to $J^{\delta,-\delta}_{3,16}$ changing basis to $\{1,b,\frac{\delta a}{\gamma}-b\}$;\\
	$(\rho \neq 0,\ \delta = 0,\ \gamma = 0):$ Isomorphic to $J^{\rho,0}_{3,16}$, no change of basis performed;\\
	$(\rho \neq 0,\ \delta = 0,\ \gamma \neq 0):$ Isomorphic to $J_{3,16}^{\rho,-\rho}$, by a change of basis to $\{1,a,a-\frac{\rho b}{\gamma}\}$;\\
	$(\rho \neq 0,\ \delta \neq 0,\ \gamma = 0):$ Isomorphic to $J^{\rho,\delta}_{3,16}$, no change of basis needed;$\\
	(\rho \neq 0,\ \delta \neq 0,\ \gamma \neq 0):$ Isomorphic to $J^{\rho-\gamma^2\delta^{-1},\delta}_{3,16}$, by a change of basis \\ to $\{1,a-\frac{\gamma b}{\delta},b\}$.
	
	Over a field of characteristic two, the following remarks will be useful.
	\begin{definition}
		An element $x$ of a unital algebra $A$ is said to have polynomial degree $n>0$ if it satisfies the equation $x^n+\rho_{n-1}x^{n-1}+\cdots+\rho_1x+\rho_01=0$, where $n$ is the smallest of such integers and $\rho_i\in K$. The degree of an element $x$ is denoted by $\partial{x}$, and the degree of an algebra is defined by $\partial{A}=\displaystyle\max_{a\in A}\{\partial a\}$.
	\end{definition}
	
	Notice that $\partial x=1$ implies $x\in \lan 1\ran$, therefore our analysis has the starting point $\partial A>1$.
	
	\begin{example}
		Let $K$ be a field and $K(\alpha)$ be a field extension of degree $n>1$, then $\alpha$ is a root of some monic polynomial $p(x)=x^n+\beta_{n-1}x^{n-1}+\cdots+\beta_0 1$ and therefore has polynomial degree $n$.
	\end{example}
	\begin{proposition}
		Let $A$ be a unital power-associative algebra of dimension $3$ that is not associative. Then $\partial{A}=2$.
	\end{proposition}
	\begin{proof}
		That $\partial A>1$ follows from dim$A>1$. Suppose $\partial A>2$ and take an element $x$ such that $\partial x>2$. The set $\{1,x,x^2\}$ is a basis for the power-associative algebra $A$, therefore it is associative.
	\end{proof}
	
	\begin{corollary}
		Any Jordan algebra of dimension $3$ that is unital and contains an element of polynomial degree $3$ is associative. 
	\end{corollary}
	
	\begin{proposition}\label{Prop5}
		Let $A$ be a unital $3$-dimensional nonassociative Jordan algebra over a field $K$ of ${\rm char}(K)=2$ with basis  $\{1,a,b\}$ satisfying $\partial a=2$ and $\partial b=2$, then $A$ satisfies one of the following conditions:
		\begin{enumerate}
			\item $A\cong J^{\alpha,\beta}_{3,16}$,
			\item $A\cong J_{3,17}$.
		\end{enumerate}
	\end{proposition}
	\begin{proof}
		Write $a^2=\alpha_1a+\alpha_01$ and $b^2=\beta_1b+\beta_01$. If $\alpha_1=0$ and $\beta_1\neq0$ or $\alpha_1\neq0$ and $\beta_1=0$ then $A=\lan 1,a+b,(a+b)^2\ran$. For $\alpha_1\neq 0$ and $\beta_1 \neq 0$ the vector $a+jb$ is cyclic for any $j\notin\{0,\alpha_1\beta_1^{-1}\}$. Therefore $\alpha_1$ and $\beta_1$ may be taken to satisfy $\alpha_1=\beta_1=0$. Since $ab=\kappa_1 1+\kappa_2a+\kappa_3b$ it follows that $(ab)^2=(\kappa_1^2+\kappa_2^2\alpha_0+\kappa_3^2\beta_0)1$. Using a change of basis to $\{1,a+\kappa_31,b+\kappa_21\}$, our basis is taken to satisfy $ab=\kappa1$, where $\kappa=\kappa_1+\kappa_2\kappa_3$. The new squares have coefficients $\alpha_0+\kappa_3^2$ and $\beta_0+\kappa_2^2$; we relabel these coefficients as $\alpha_0$ and $\beta_0$. Investigating the conditions over $\kappa,\alpha_0,\beta_0$ leads to the following.\\

		Taking $\alpha = \alpha_0$ and $\beta= \beta_0$, the product of the algebra is codified in the triple
		$(\alpha, \beta, \kappa)$. Hence, if $\kappa=0$ and $(\alpha, \beta)\neq (0,0)$, then $A$ is the algebra $J^{\alpha,\beta}_{3,16}$. The triple $(0,0,0)$ produces an associative algebra. Moreover, if $\kappa \ne 0 $ and
		$( \alpha, \beta ) \neq (0,0)$, then under a suitable change of basis we obtain an algebra of type $J^{\alpha^\prime,\beta^\prime}_{3,16}$. \\
		
		In fact, for the triples $(0,\beta \neq 0, \kappa \neq 0)$ and 
		$(\alpha \neq 0,\beta, \kappa\neq 0)$, take the change of basis $\{1,\beta^{-1}\kappa b+a,b\}$ and 
		$\{1, a,b+\alpha^{-1}\kappa a\}$ respectively to obtain the algebras  $J^{\kappa^2\beta^{-1},\beta}_{3,16}$ and $J_{3,16}^{\alpha,\beta+\kappa^2\alpha^{-1}}$. Finally, for the triple $(0,0, \kappa\neq 0)$ we obtain the algebra $J_{3,17}$.
	\end{proof}
	
	\subsection{Isomorphism classes over a field of characteristic not two}
	\leavevmode\par
	\begin{proposition} Algebras $J^{\alpha,\beta}_{3,16}$ and $J^{\chi,\lambda}_{3,16}$ are isomorphic if and only if there are $\kappa_1,\kappa_2,\kappa_3,\kappa_4\in K$ with $\kappa_1\kappa_4-\kappa_3\kappa_2\neq 0$, $\kappa_1\kappa_3\chi+\kappa_2\kappa_4\lambda=0$ and 
		$$\begin{pmatrix}
			\alpha \\
			\beta
		\end{pmatrix}= 
		\begin{pmatrix}
			\kappa_1^2 & \kappa_2^2 \\
			\kappa_3^2 & \kappa_4^2
		\end{pmatrix} 
		\begin{pmatrix}
			\chi \\
			\lambda
		\end{pmatrix}.$$ 
	\end{proposition}
	\begin{proof}
		Suppose there is an isomorphism:
		\begin{center}
			$ h\colon J^{\alpha,\beta}_{3,16}=\lan 1,a,b \mid a^2 =\alpha 1,b^2=\beta 1\ran \to  J^{\chi,\lambda}_{3,16}= \lan 1,x,y \mid x^2=\chi 1,y^2=\lambda 1\ran 
			$
		\end{center}
		The homomorphism $h$ must satisfy:
		\begin{enumerate}
			\item $h(1)=1$, $h(a)=c_11+c_2x+c_3y$, and  $h(b)=c_41+c_5x+c_6y;$
			\item $\det ([h]) = \left|\begin{array}{ccc}
				1 & c_1 & c_4\\
				0 & c_2 & c_5\\
				0 & c_3 & c_6
			\end{array}\right|=c_2c_6-c_3c_5\neq 0 $;
			\item $h(a)^2=\alpha1$, $h(b)^2=\beta 1$, $h(ab)=0$.
		\end{enumerate}
		From where we deduce that:
		$$
		\begin{cases}
			c_2c_6-c_5c_3\neq 0; \\
			c_1c_2=0;\\
			c_1c_3=0;\\
			c_1^2+c_2^2\chi+c_3^2\lambda=\alpha; \\
			c_4c_5=0; \\
			c_4c_6=0;\\
			c_4^2+c_5^2\chi+c_6^2\lambda = \beta; \\
			c_1c_4+c_2c_5\chi+c_3c_6\lambda=0.
		\end{cases}
		$$
		Since $c_1\neq 0$ implies $h(a)=c_1h(1)$ and a similar result holds for $c_4\neq 0$ we conclude that $c_1=c_4=0$ and the above conditions rapidly imply the result. Conversely, if such $\kappa_1,\ldots,\kappa_4 \in K$ exist, a change of basis to $\{1,\kappa_1 x+\kappa_2y,\kappa_3x+\kappa_4 y\}$ proves the isomorphism.
	\end{proof}
	
	\begin{proposition} Algebras $J_{3,16}^{\alpha,\beta} $ and $J_{3,17}$ are isomorphic if and only if $-\alpha\beta=r^2\neq 0$ for some $r\in K$.
	\end{proposition}
	\begin{proof}
		Suppose there is an isomorphism
		$$
		h:J_{3,16}^{\alpha,\beta}
		=
		\left\langle 1,a,b \mid
		a^2=\alpha 1,\;
		b^2=\beta 1
		\right\rangle\to
		J_{3,17}
		=
		\left\langle 1,a',b'\mid 
		a'b'=1
		\right\rangle.
		$$
		then it must satisfy:
		
		\begin{enumerate}
			\item$h(1)=1.$
			\item
			$h(a)=\alpha_1+\alpha_2a'+\alpha_3b'$, $
			h(b)=\beta_1+\beta_2a'+\beta_3b'$.
			\item Since $a^2=\alpha 1$, we have $h(a)^2=h(a^2)=\alpha h(1)=\alpha$.
			Hence
			$$
			\alpha=
			\left(\alpha_1^2+2\alpha_2\alpha_3\right)
			+
			\left(2\alpha_1\alpha_2\right)a'
			+
			\left(2\alpha_1\alpha_3\right)b'.
			$$
			\item Since $h$ is an isomorphism,
			$
			\det
			\begin{pmatrix}
				1 & \alpha_1 & \beta_1\\
				0 & \alpha_2 & \beta_2\\
				0 & \alpha_3 & \beta_3 
			\end{pmatrix}
			\neq 0.
			$
			
			\item
			Since $b^2=\beta 1$, we have
			$h(b)^2=\beta$,
			and therefore
			$$
			\beta
			=
			\left(\beta_1^2+2\beta_2\beta_3\right)
			+
			\left(2\beta_1\beta_2\right)a'
			+
			\left(2\beta_1\beta_3\right)b'.
			$$
			
			\item
			Since $ab=0$ in
			$J_{3,16}^{\alpha,\beta}$, $h(ab)=0$.
			Thus
			$$
			0
			=
			\left(
			\alpha_1\beta_1
			+
			\alpha_2\beta_3
			+
			\alpha_3\beta_2
			\right)
			+
			\left(
			\alpha_1\beta_2
			+
			\alpha_2\beta_1
			\right)a'
			+
			\left(
			\alpha_1\beta_3
			+
			\alpha_3\beta_1
			\right)b'.
			$$
		\end{enumerate}
		Consequently, condition (3) gives
		$$
		\begin{cases}
			\alpha=\alpha_1^2+2\alpha_2\alpha_3,\\
			2\alpha_1\alpha_2=0,\\
			2\alpha_1\alpha_3=0.
		\end{cases}
		$$
		Condition (5) gives
		$$
		\begin{cases}
			\beta=\beta_1^2+2\beta_2\beta_3,\\
			0=2\beta_1\beta_2,\\
			0=2\beta_1\beta_3.
		\end{cases}
		$$
		Condition $(6)$ gives
		$$
		\begin{cases}
			\alpha_1\beta_1+\alpha_2\beta_3+\alpha_3\beta_2=0,\\
			\alpha_1\beta_2+\alpha_2\beta_1=0,\\
			\alpha_1\beta_3+\alpha_3\beta_1=0.
		\end{cases}
		$$
		If $\alpha_1\neq 0$, then, since ${\rm char(K)}\neq 2$, 
		$\alpha_2=\alpha_3=0$, and therefore
		$\{1,h(a)\}$ would be linearly dependent.
		Hence $\alpha_1=0$. By the same argument, $\beta_1=0$. Thus the previous conditions reduce to
		$$
		\begin{cases}
			\alpha=2\alpha_2\alpha_3, \\
			\beta=2\beta_2\beta_3,\\
			\alpha_2\beta_3+\alpha_3\beta_2=0,
		\end{cases}
		$$
		and condition (4) becomes $\alpha_2\beta_3-\beta_2\alpha_3\neq 0.$\\
		
		Suppose $\alpha=0$. Then either $\alpha_2=0$ or $\alpha_3=0$. If
		$
		\alpha_2=0,
		$
		then,
		$\alpha_3\beta_2=0$. Since $\alpha_1=0$, if $\alpha_3=0$, then
		$h(a)=0$, a contradiction. Hence $\beta_2=0.$ But then
		$\alpha_2\beta_3-\beta_2\alpha_3=0$, contradicting condition (4). In a similar manner, $\alpha_3=0$ is also impossible. Therefore $\alpha\neq 0$.
		Similarly, by the same argument, $\beta\neq 0$. Since
		$\alpha\neq 0,$ we have $\alpha_2\neq 0,\ \alpha_3\neq 0$, by the same reasoning $\beta_2\neq 0, \beta_3\neq 0$. Using that $\alpha_2\beta_3=-\alpha_3\beta_2$. We conclude
		\begin{align*}
			\alpha\beta
			&=
			4\alpha_2\alpha_3\beta_2\beta_3=
			4(\alpha_2\beta_3)(\beta_2\alpha_3)=
			4(\alpha_2\beta_3)(-\alpha_2\beta_3).
		\end{align*}
		Hence
		$-\alpha\beta=2^2\alpha_2^2\beta_3^2$, and consequently
		$-\alpha\beta=\left(2\alpha_2\beta_3\right)^2
		$. Since the scalars involved are nonzero, $-\alpha\beta$ is a nonzero square in $K$. For the converse, a change of basis to $\{1,a+\frac{r}{\beta}b,\frac{a}{2\alpha}
		-\frac{r}{2\alpha\beta}b\}$ applied to $J_{3,16}^{\alpha,\beta}$ proves the isomorphism.
	\end{proof}
	
	\subsection{Isomorphism classes over a field of characteristic two}
	\leavevmode\par 	
	\begin{proposition} Algebras $J^{\alpha,\beta}_{3,16},J^{\chi,\lambda}_{3,16}$ are isomorphic if and only if there are $\kappa_1,\ldots,\kappa_6\in K$ such that $\kappa_1\kappa_4+\kappa_2\kappa_5\chi+\kappa_3\kappa_6\lambda=0,\kappa_1\kappa_5+\kappa_2\kappa_4=0,\kappa_3\kappa_4+\kappa_1\kappa_6=0,\kappa_2\kappa_6+\kappa_3\kappa_5\neq 0,\kappa_1^2+\kappa_2^2\chi+\kappa_3^2\lambda=\alpha, \kappa_4^2+\kappa_5^2\chi+\kappa_6^2\lambda=\beta$.
	\end{proposition}
	\begin{proof}
		Suppose there is an isomorphism:
		\begin{center}
			$ h\colon J^{\alpha,\beta}_{3,16}=\lan 1,a,b \mid a^2 =\alpha 1,b^2=\beta 1\ran\to \lan 1,x,y \mid x^2=\chi 1,y^2=\lambda 1\ran =J^{\chi,\lambda}_{3,16}$
		\end{center}
		The homomorphism $h$ must satisfy:
		\begin{enumerate}
			\item $h(1)=1$, $h(a)=\kappa_1+\kappa_2x+\kappa_3y$,   $h(b)=\kappa_41+\kappa_5x+\kappa_6y$;
			\item $\det ([h]) = \left|\begin{array}{ccc}
				1 & \kappa_1 & \kappa_4\\
				0 & \kappa_2 & \kappa_5\\
				0 & \kappa_3 & \kappa_6
			\end{array}\right|=\kappa_2\kappa_6-\kappa_3\kappa_5\neq 0 $;
			\item $h(a)^2=\alpha1$, $h(b)^2=\beta 1$, $h(ab)=0$.
		\end{enumerate}
		From where we deduce that:
		$$
		\begin{cases}
			\kappa_2\kappa_6-\kappa_5\kappa_3\neq 0; \\
			\kappa_1^2+\kappa_2^2\chi+\kappa_3^2\lambda=\alpha; \\
			\kappa_4^2+\kappa_5^2\chi+\kappa_6^2\lambda = \beta; \\
			\kappa_1\kappa_4+\kappa_2\kappa_5\chi+\kappa_3\kappa_6\lambda=0;\\
			\kappa_1\kappa_5+\kappa_2\kappa_4=0;\\
			\kappa_1\kappa_6+\kappa_3\kappa_4=0.
		\end{cases}
		$$
		This proves the implication. The converse is proved considering a change of basis to $\{1,\kappa_11+\kappa_2x+\kappa_3y,\kappa_41+\kappa_5x+\kappa_6y\}$.
	\end{proof}
	
	\begin{proposition} Algebras $J^{\alpha,\beta}_{3,16}$ and $J_{3,17}$ are not isomorphic.
	\end{proposition}
	\begin{proof}
		An isomorphism would give elements
		$$
		u=\mu_11+\mu_2a+\mu_3b,\quad v=\nu_11+\nu_2a+\nu_3b
		$$
		such that $\{1,u,v\}$ is a basis, $u^2=v^2=0$ and $uv=1$. Linear independence is ruled by the equation $\mu_2\nu_3+\mu_3\nu_2\neq 0$. The coefficients of $a$ and $b$ in the product $uv$ are $\mu_1\nu_2+\nu_1\mu_2=0$ and $\mu_1\nu_3+\mu_3\nu_1=0$. The linear independence relation forces $\mu_1=
		\nu_1=0$. The square condition provides:
		$$\begin{pmatrix}
			0 \\
			0
		\end{pmatrix}= 
		\begin{pmatrix}
			\mu_2^2 & \mu_3^2 \\
			\nu_2^2 & \nu_3^2
		\end{pmatrix} 
		\begin{pmatrix}
			\alpha \\
			\beta
		\end{pmatrix}.$$ 
		Since the determinant of the first matrix of the right-hand side is nonzero, we must have $\alpha=\beta=0$. This implies that $uv=0$, contradicting $uv=1$.
	\end{proof}
	\subsection{Approach for the nonunital case over a field of characteristic not two.}
	
	\begin{proposition}\label{proposition5} Let $A$ be a $3$-dimensional nonunital Jordan algebra over a field $K$ with $char K\neq2$, and suppose that $A$ is not a nilalgebra. Then the following hold: 
		\begin{enumerate} 
			\item[(a)]  $A$ has type $(1,2,0)$ if and only if $$ A\cong J_{3,18} =\lan e,a,b\mid e^2=e,\ ea=\frac{a}{2},\ eb=\frac{b}{2}\ran. $$ 
			
			\item[(b)] $A$ has type $(1,1,1)$ if and only if it is isomorphic to one of the following algebras: $$ J_{3,19} =\lan e,a,a^2\mid e^2=e,\ ea=\frac{a}{2}\ran,  \  J_{3,20} =\lan e,a,b\mid e^2=e,\ ea=\frac{a}{2}\ran. $$ 
			
			\item[(c)] $A$ has type $(2,1,0)$ if and only if it is isomorphic to one of the following algebras: $$ J_{3,21} =\lan e,b,b^2\mid e^2=e,\ eb=\frac{b}{2},\ eb^2=b^2\ran, $$ \vspace{-0,5cm} $$ J_{3,22} =\lan e,a,b\mid e^2=e,\ ea=a,\ eb=\frac{b}{2}\ran, $$ \vspace{-0,5cm} $$ J_{3,23} =\lan e,a,b\mid e^2=e,\ ea=a,\ eb=\frac{b}{2},\ a^2=e,\ ab=\frac{b}{2}\ran. $$
		\end{enumerate} 
	\end{proposition}
	\begin{proof}
		Recall that a linearization of the Jordan identity gives
		$$
		(wz,y,x)+(xz,y,w)+(xw,y,z)=0.
		$$
		We recall that throughout the proof, the idempotent $e$ is chosen so that $\dim A_1$ is maximal among all Peirce decompositions of $A$ and, among those attaining the same value of $\dim A_1$, $\dim A_{1/2}$ is maximal.\\
		
		$(a)\Rightarrow$ Let $A_1=\lan e\ran$ and $A_{1/2}=\lan a,b\ran$. By the Peirce multiplication rules, the products are
		$$
		e^2=e, \  ea=\frac{a}{2}, \  eb=\frac{b}{2}, \ 
		a^2=\kappa_1e, \  b^2=\kappa_2e, \  ab=\kappa_3e.
		$$
		From $(a^2,a,a)=0$ and $(b^2,b,b)=0$ we obtain $\kappa_1=\kappa_2=0$. Hence $(a+b)^2=2\kappa_3e$, and the identity $((a+b)^2,e,a+b)=0$ gives $\kappa_3=0$. Therefore the only nonzero products are
		$$
		e^2=e, \  ea=\frac{a}{2}, \  eb=\frac{b}{2},
		$$
		and consequently $A\cong J_{3,18}$.
		
		Conversely, a direct computation shows that $J_{3,18}$ satisfies the Jordan identity. Moreover, every nonzero idempotent is of the form $f=e+\alpha a+\beta b$, and its Peirce decomposition has dimensions $(1,2,0)$. Thus the maximal type of $J_{3,18}$ is $(1,2,0)$.\\

		$(b)\Rightarrow$ Since $A$ has type $(1,1,1)$, let $A_1=\lan e\ran$, $A_{1/2}=\lan a\ran$ and $A_0=\lan b\ran$. Then
		$$
		e^2=e, \  ea=\frac{a}{2}, \  eb=0, \ 
		b^2=\kappa b, \  a^2=\kappa_1e+\kappa_2b, \  ab=\kappa_3a.
		$$
		If $\kappa\neq0$, after replacing $b$ by $\kappa^{-1}b$ we may assume that $b^2=b$. In this case $e+b$ is an idempotent and both $e$ and $b$ belong to $A_1(e+b)$, since $(e+b)e=e$ and $(e+b)b=b$. It follows that $\dim A_1(e+b)\geq2$, contradicting the maximality of $\dim A_1(e)=1$. Hence $\kappa=0$. \\
		
		Now, from $(a^2,e,a)=0$ and $(ab,a,e)+(eb,a,a)+(ea,a,b)=0$
		we obtain, respectively, $\kappa_1-2\kappa_2\kappa_3=0$
		and $\kappa_2\kappa_3=0$.
		Therefore $\kappa_1=0$. If $\kappa_2=0$ and $\kappa_3\neq0$, then $a^2=b^2=0$ and $ab=\kappa_3a$, so that $((a+b)^2,b,a+b)=2\kappa_3^3a\neq0$, contradicting the Jordan identity. Thus this case cannot occur.
		
		If $\kappa_2\neq0$, necessarily $\kappa_3=0$, and $a^2=\kappa_2b\neq0$. Hence $\{e,a,a^2\}$ is a basis of $A$, and with respect to this basis we obtain $A\cong J_{3,19}$. If $\kappa_2=\kappa_3=0$, the only nonzero products are $e^2=e$ and $ea=a/2$, and therefore $A\cong J_{3,20}$.
		
		Conversely, both $J_{3,19}$ and $J_{3,20}$ satisfy the Jordan identity. In $J_{3,19}$, every nonzero idempotent is of the form $f=e+\alpha a+\alpha^2a^2$, and its Peirce decomposition has dimensions $(1,1,1)$. In $J_{3,20}$, every nonzero idempotent is of the form $f=e+\alpha a$, and the corresponding Peirce decomposition again has dimensions $(1,1,1)$. Hence both algebras have maximal type $(1,1,1)$.\\
		
		$(c)\Rightarrow$ Suppose that $A$ has type $(2,1,0)$ and fix a basis such that $A_1=\lan e,a\ran$ and $A_{1/2}=\lan b\ran$. The multiplication has the form
		$$
		e^2=e, \  ea=a, \  eb=\frac{b}{2}, \ 
		a^2=\kappa_1e+\kappa_2a, \  ab=\kappa b, \ 
		b^2=\kappa_3e+\kappa_4a.
		$$
		Replacing $a$ by $a-\frac{\kappa_2}{2}e$ and relabelling the constants, we may assume that $a^2=\kappa_1e$. Thus we may write
		$$
		a^2=\kappa_1e, \  ab=\kappa b, \ 
		b^2=\kappa_3e+\kappa_4a.
		$$
		
		Applying the linearized Jordan identity to
		$$
		(ea,a,b)+(ba,a,e)+(be,a,a)=0
		$$
		gives $\kappa_1=4\kappa^2$. On the other hand, $(b^2,e,b)=0$ gives
		$\kappa_3+2\kappa\kappa_4=0$, while
		$$
		(ea,b,b)+(ba,b,e)+(be,b,a)=0
		$$
		gives $\kappa_3-2\kappa\kappa_4=0$. Since ${\rm char}(K)\neq2$, these two relations imply $\kappa_3=0$ and $\kappa\kappa_4=0$. If $\kappa=0$, then $\kappa_1=0$, and consequently
		$$
		a^2=0, \  ab=0, \  b^2=\kappa_4a.
		$$
		When $\kappa_4=0$ we obtain $J_{3,22}$. If $\kappa_4\neq0$, then $b^2=\kappa_4a\neq0$ and $\{e,b,b^2\}$ is a basis of $A$. Since $eb=b/2$ and $e b^2=b^2$, this gives $A\cong J_{3,21}$.\\
		
		Suppose now that $\kappa\neq0$. Then $\kappa_4=0$, so $b^2=0$, and $\kappa_1=(2\kappa)^2$. Replacing $a$ by $a'=(2\kappa)^{-1}a$, we obtain $(a')^2=e$ and $a'b=b/2$. Thus, with respect to the basis $\{e,a',b\}$:
		$$
		e^2=e, \  ea'=a', \  eb=\frac{b}{2}, \ 
		(a')^2=e, \  a'b=\frac{b}{2},
		$$
		and hence $A\cong J_{3,23}$.
		
		Conversely, a direct verification shows that $J_{3,21}$, $J_{3,22}$ and $J_{3,23}$ satisfy the Jordan identity. In each case the idempotent $e$ satisfies $\dim A_1(e)=2$ and $\dim A_{1/2}(e)=1$. Since these algebras are nonunital, no idempotent $f$ can satisfy $\dim A_1(f)=3$, for otherwise $L_f=Id_A$ and $f$ would be an identity element. Hence $\dim A_1=2$ is maximal and the maximal type of each of these algebras is $(2,1,0)$.
	\end{proof}
	
	Algebras of type $(2,0,1)$ are direct sums of associative algebras. Indeed,
	$A_1$ is a two-dimensional unital Jordan algebra and hence is associative,
	while $A_0$ is one-dimensional and $A_1A_0=0$. Moreover, since $A$ is
	nonunital, $A_0$ cannot contain a nonzero idempotent.\\
	
	For type $(1,0,2)$, we have $A=A_1\oplus A_0,A_1=\lan e\ran$. By the maximality of $\dim A_1$, the algebra $A_0$ cannot contain a nonzero
	idempotent. Indeed, if $f\in A_0$ were a nonzero idempotent, then $(e+f)^2=e+f$ and $\lan e,f\ran\subseteq A_1(e+f)$, contradicting the maximality of $\dim A_1(e)=1$. In particular, the case
	$A_0\cong
	\lan f,a\mid f^2=f,\;fa=\frac{a}{2}\ran$ does not have maximal type $(1,0,2)$. In fact, considering the basis $\{e+f,f-e,a\}$, we obtain
	$A\cong J_{3,23}$ whose maximal type is $(2,1,0)$. This shows that $A$ is a direct sum of orthogonal associative algebras, hence associative.
	
	\subsection{Approach for the nonunital case over characteristic two.}
	$$J_{3,24}=\langle a,b,c\ |\ a^2=a,bc=b\rangle$$
	Given an idempotent $e$ there are only types $(2,1),(1,2)$ to be considered. Since $A_1$ is unital, the type $(2,1)$ always leads to an associative algebra. The only two-dimensional Jordan algebra that is not associative is $J_{2,7}$, therefore the unique algebra that is not associative and has type $(1,2)$ is given by $J_{1,2}\oplus J_{2,7}\cong J_{3,24}$.
	
	\section{Geometric classification over algebraically closed fields of characteristic two}
	Let $V$ be an $n$-dimensional $K$-vector space with a fixed basis
	$\{e_1,e_2,\ldots,e_n\}$. To endow $V$ with a Jordan $K$-algebra
	structure, $(J,\cdot)$, it suffices to specify $n^3$ structure
	constants $\alpha_{ij}^k\in K$, namely,
	$e_i\cdot e_j=\sum_{k=1}^{n}\alpha_{ij}^k e_k$,
	$i,j\in\{1,2,\ldots,n\}$. If a point $(\alpha_{ij}^k)\in K^{n^3}$
	defines a Jordan algebra its coordinates satisfy certain polynomial equations
	which follow from commutativity and the Jordan identity. These equations define
	an algebraic variety $\mathrm{Jor}_n$ in $K^{n^3}$. The points
	$(\alpha_{ij}^k)\in\mathrm{Jor}_n$ are in one-to-one correspondence with
	$n$-dimensional Jordan $K$-algebra structures on $V$ in the fixed basis.
	A change of basis in $\mathcal{J}$ is given by the general linear group
	$G=\mathrm{GL}(V)$ action on $\mathrm{Jor}_n$, namely:
	\begin{equation}
		g(\mathcal{J},\cdot)\mapsto(\mathcal{J},\cdot_g),\qquad
		x\cdot_g y=g\bigl(g^{-1}x\cdot g^{-1}y\bigr),
		\tag{2}
	\end{equation}
	for any $\mathcal{J}\in\mathrm{Jor}_n$, $g\in G$ and $x,y\in V$.
	The $G$-orbit of a Jordan algebra $\mathcal{J}$, that is, the set of all
	images of $\mathcal{J}$ under the action of $G$, is denoted by
	$\mathcal{J}^{G}$. The set of different $G$-orbits of this action is in
	one-to-one correspondence with the set of isomorphism classes of
	$n$-dimensional Jordan algebras. The Zariski topology of the affine space
	provides the following relation of orbits, namely, we say that $\mathcal{J}_1$
	is a deformation of $\mathcal{J}_2$, if the orbit $\mathcal{J}_2^{G}$ is contained in the Zariski
	closure of the orbit $\mathcal{J}_1^{G}$. We denote it by
	$\mathcal{J}_1\to\mathcal{J}_2$. A Jordan algebra $\mathcal{J}$ is called
	rigid if its $G$-orbit is a Zariski-open set in $\mathrm{Jor}_n$. Rigid
	algebra $\mathcal{J}$ does not admit non-trivial deformations, indeed if
	$\mathcal{J}_1\to\mathcal{J}$ then
	$\overline{\mathcal{J}_1^{G}}\cap\mathcal{J}^{G}\neq\emptyset$
	and therefore $\mathcal{J}_1\cong\mathcal{J}$.\\
	
	The main goal of this section is to determine the rigid algebras in dimension 2 and 3 as well as determine the irreducible components of the varieties $\rm{Jor}_2$ and $\rm{Jor}_3$ over an algebraically closed field of characteristic two.\\
	
	As a criterion to determine the non-existence of deformations we establish the notion of Characteristic and Alicyclic algebras. Before we introduce these definitions we prove some auxiliary results and make a few remarks. Let $A\in M_n(K)$. Recall that the characteristic polynomial of $A$ can be written as
	$$
	\chi_A(t)
	=
	t^n-E_1(A)t^{n-1}+\cdots
	+(-1)^{n-1}E_{n-1}(A)t+(-1)^nE_n(A),
	$$
	where $E_k(A)$ denotes the sum of all principal minors of order $k$ of $A$; see, for instance, \cite[p.53, Eq.(1.2.13)]{MA}.
	
	\begin{proposition}
		Let $A\in M_n(K)$. For each $1\leq k\leq n$, the coefficient
		$E_k(A)$
		of the characteristic polynomial of $A$ is a homogeneous polynomial of degree $k$ in the entries of $A$.
	\end{proposition}
	
	\begin{proof}
		By the permutation expansion of the determinant,
		$$
		\det(A)
		=\sum_{\sigma\in S_n}
		\operatorname{sgn}(\sigma)
		\prod_{i=1}^n a_{i,\sigma(i)},
		$$
		see \cite[p.~9, Eq.~(0.3.2.1)]{MA}. Each summand is a monomial of degree $n$ in the entries of $A$. Hence $\det(A)$ is a homogeneous polynomial of degree $n$ in the entries of $A$.
		
		Now, every principal minor of order $k$ of $A$ is the determinant of a $k\times k$ principal submatrix of $A$. By the same argument, each such determinant is a homogeneous polynomial of degree $k$ in the entries of $A$. Since $E_k(A)$ is the sum of all principal minors of order $k$, it is also a homogeneous polynomial of degree $k$. Therefore $E_k(A)$ is a homogeneous polynomial of degree $k$ in the entries of $A$.
	\end{proof}
	We now apply this observation to multiplication operators of a finite-dimensional algebra. Let $A$ be an $n$-dimensional algebra over $K$, and fix a basis $\mathcal B=\{e_1,\ldots,e_n\}$. Let $a=x_1e_1+\cdots+x_ne_n\in A$, and let $\alpha_{ij}^k\in K$ be the structure constants of $A$ with respect to $\mathcal B$, so that
	$e_ie_j=\sum_{k=1}^n\alpha_{ij}^ke_k$. Consider the left multiplication operator $L_a:A\longrightarrow A$, $L_a(v)=av$. For every basis element $e_j$, we have
	\begin{align*}
		L_a(e_j)=ae_j =\left(\sum_{i=1}^n x_i e_i\right)e_j=\sum_{i=1}^n x_i(e_ie_j)=\sum_{k=1}^n
		\left(
		\sum_{i=1}^n x_i\alpha_{ij}^k
		\right)e_k.
	\end{align*}
	Thus, if
	$$
	B=[L_a]_{\mathcal B}^{\mathcal B}=(b_{\ell j})_{\ell,j=1}^n,
	$$
	then
	$$
	b_{\ell j}
	=\sum_{i=1}^n x_i\alpha_{ij}^\ell.
	$$
	Consequently, every entry of the matrix $B=[L_a]_{\mathcal B}^{\mathcal B}$ is a homogeneous polynomial of degree $1$ in the coordinates $x_1,\ldots,x_n$ of $a$.
	
	\begin{proposition}\label{prop111}
		Let
		$B=[L_a]_{\mathcal B}^{\mathcal B}$.
		For every $1\leq k\leq n$, the coefficient of $t^{n-k}$ in the characteristic polynomial of $B$ is a homogeneous polynomial of degree $k$ in $x_1,\ldots,x_n$.
	\end{proposition}
	\begin{proof}
		The coefficient of $t^{n-k}$ in $\chi_B(t)$ is $(-1)^kE_k(B)$.
		By the previous proposition, $E_k(B)$ is homogeneous of degree $k$ in the entries of $B$. Since every entry
		$$
		b_{\ell j}=\sum_{i=1}^n x_i\alpha_{ij}^\ell,
		$$
		is homogeneous of degree $1$ in $x_1,\ldots,x_n$, it follows that $E_k(B)$, and therefore $(-1)^kE_k(B)$, is homogeneous of degree $k$ in $x_1,\ldots,x_n$.
	\end{proof}
	We may therefore write
	$$
	\chi_{L_a}(t)
	=
	t^n
	+s_1(x_1,\ldots,x_n)t^{n-1}
	+\cdots
	+s_{n-1}(x_1,\ldots,x_n)t
	+s_n(x_1,\ldots,x_n),
	$$
	where, for each $k$,
	$s_k(x_1,\ldots,x_n)=(-1)^kE_k(B)$
	is homogeneous of degree $k$ in $x_1,\ldots,x_n$. We examine the dependence of the coefficients of these polynomials on the structure constants of $A$.
	
	\begin{proposition}\label{teohomopoly}
		For every $1\leq k\leq n$, considering
		$s_k(x_1,\ldots,x_n)$ as a polynomial in $x_1,\ldots,x_n$, every coefficient of $s_k$ is a homogeneous polynomial of degree $k$ in the structure constants
		$\alpha_{ij}^\ell$ of the algebra $A$.
	\end{proposition}
	\begin{proof}
		Since
		$$
		s_k(x_1,\ldots,x_n)=(-1)^kE_k(B),
		$$
		it is enough to consider an arbitrary principal minor of order $k$ of $B$.
		
		Let
		$$
		1\leq i_1<\cdots<i_k\leq n,
		$$
		and denote by
		$$
		B[i_1,\ldots,i_k\mid i_1,\ldots,i_k]
		$$
		the principal submatrix of $B$ whose rows and columns are indexed by
		$i_1,\ldots,i_k$. Its determinant is
		$$
		\det B[i_1,\ldots,i_k\mid i_1,\ldots,i_k]
		=
		\sum_{\sigma\in S_k}
		\operatorname{sgn}(\sigma)
		\prod_{j=1}^k
		b_{i_j,i_{\sigma(j)}}.
		$$
		Using
		$$
		b_{\ell r}
		=
		\sum_{p=1}^n x_p\alpha_{pr}^\ell,
		$$
		we obtain
		\begin{align*}
			\det B[i_1,\ldots,i_k\mid i_1,\ldots,i_k]
			&=
			\sum_{\sigma\in S_k}
			\operatorname{sgn}(\sigma)
			\prod_{j=1}^k
			\left(
			\sum_{p=1}^n
			x_p\alpha_{p\,i_{\sigma(j)}}^{i_j}
			\right)\\
			&=
			\sum_{(p_1,\ldots,p_k)\in\{1,\ldots,n\}^k}
			\left(
			\sum_{\sigma\in S_k}
			\operatorname{sgn}(\sigma)
			\prod_{j=1}^k
			\alpha_{p_j\,i_{\sigma(j)}}^{i_j}
			\right)
			x_{p_1}\cdots x_{p_k}.
		\end{align*}
		
		For each fixed $(p_1,\ldots,p_k)$, the expression
		$$
		\sum_{\sigma\in S_k}
		\operatorname{sgn}(\sigma)
		\prod_{j=1}^k
		\alpha_{p_j\,i_{\sigma(j)}}^{i_j}
		$$
		is a homogeneous polynomial of degree $k$ in the structure constants of $A$, since every product occurring in the sum contains exactly $k$ structure constants.
		
		There may, of course, be different tuples $(p_1,\ldots,p_k)$ giving rise to the same monomial in the commuting variables
		$x_1,\ldots,x_n$. After collecting equal monomials, their coefficient is therefore a sum of homogeneous polynomials of degree $k$ in the structure constants. Hence it is again homogeneous of degree $k$ in the structure constants.
		
		Thus every coefficient, with respect to the variables $x_1,\ldots,x_n$, of every principal minor of order $k$ of $B$ is a homogeneous polynomial of degree $k$ in the structure constants of $A$. Since $E_k(B)$ is the sum of all such principal minors, the same conclusion holds for
		$$
		s_k(x_1,\ldots,x_n)=(-1)^kE_k(B).
		$$
		Therefore every coefficient of $s_k$ is a homogeneous polynomial of degree $k$ in the structure constants $\alpha_{ij}^\ell$.
	\end{proof}
	
	\begin{definition}
		We say that a nonassociative algebra $A$ with basis $\{e_1, \dots, e_n\},n>1$, is $alicyclic$ if it satisfies:
		$$
		\forall a\in A\big( \det(a,\dots,a^n)=0\big),
		$$
		where the powers of $a$ are defined by $a^1=a,a^k=aa^{k-1}$ for $k\geq 2$ and by $\det(a,\dots,a^n)=0$ we mean the determinant of the matrix whose entries are the coordinates of the vectors $a,\ldots,a^n$.
	\end{definition}
	
	\begin{proposition}\label{propositionAlicylic}
		Let $A$ and $B$ be two algebras. If $A$ is alicyclic and there is a deformation from $A$ to $B$, then $B$ is alicyclic.
	\end{proposition}
	\begin{proof}
		Let $A$ be an $n$-dimensional nonassociative algebra over a field $K$ and fix a basis
		$
		\mathcal B=\{e_1,\ldots,e_n\}
		$
		of $A$. We write the structure constants of $A$ with respect to this basis as
		$$
		e_i e_j=\sum_{k=1}^n \alpha_{ij}^k e_k,
		\ 1\leq i,j\leq n.
		$$
		Here the lower indices $i,j$ indicate the two basis elements being multiplied, while the upper index $k$ indicates the coefficient of the basis vector $e_k$ in the product.
		
		Let
		$$
		x=\sum_{i=1}^n x_i e_i\in A.
		$$
		Associated with $x$, we consider the left multiplication operator
		$$
		L_x:A\to A,\ L_x(a)=xa.
		$$
		Since the product of $A$ is bilinear, $L_x$ is a $K$-linear map.
		
		We first compute the matrix of $L_x$ in the fixed basis $\mathcal B$. For a basis element $e_j$, we have
		$$
		L_x(e_j)=xe_j
		=
		\left(\sum_{i=1}^n x_i e_i\right)e_j
		=
		\sum_{i=1}^n x_i(e_i e_j).
		$$
		Using the structure constants,
		$$
		e_i e_j=\sum_{k=1}^n \alpha_{ij}^k e_k,
		$$
		we obtain
		$$
		L_x(e_j)
		=
		\sum_{i=1}^n x_i
		\left(
		\sum_{k=1}^n \alpha_{ij}^k e_k
		\right)
		=
		\sum_{k=1}^n
		\left(
		\sum_{i=1}^n \alpha_{ij}^k x_i
		\right)e_k.
		$$
		Therefore, if $[L_x]_{\mathcal B}$ denotes the matrix of $L_x$ with respect to the basis $\mathcal B$, then its entry in row $k$ and column $j$ is
		$$
		([L_x]_{\mathcal B})_{kj}
		=
		\sum_{i=1}^n \alpha_{ij}^k x_i.
		$$
		The column index $j$ corresponds to the basis vector $e_j$ to which $L_x$ is applied, while the row index $k$ records the coefficient of $e_k$ in the result.
		
		Now define
		$$
		v_0=x,
		\
		v_m=L_x^m(x)
		\ \text{for }m\geq 1.
		$$
		Thus
		$$
		v_1=L_x(x)=xx,
		$$
		$$
		v_2=L_x^2(x)=x(xx),
		$$
		$$
		v_3=L_x^3(x)=x(x(xx)),
		$$
		and, in general, $v_m$ is the right-nested product obtained by applying left multiplication by $x$ exactly $m$ times to $x$. Since the algebra is not assumed to be associative, this parenthesization is part of the definition.\\
		
		For each $m\geq 0$, write
		$$
		v_m=\sum_{k=1}^n v_m^k e_k.
		$$
		We shall show that $v_m^k$ is a homogeneous polynomial of degree $m+1$ in the variables $x_1,\ldots,x_n$, whose coefficients are polynomials in the structure constants $\alpha_{ij}^k$.
		
		For $m=0$, we have
		$$
		v_0=x=\sum_{k=1}^n x_k e_k,
		$$
		so
		$$
		v_0^k=x_k.
		$$
		Thus $v_0^k$ is homogeneous of degree $1$ in the variables $x_1,\ldots,x_n$.
		
		For $m=1$, we compute
		$$
		v_1=xx
		=
		\left(\sum_{i_0=1}^n x_{i_0}e_{i_0}\right)
		\left(\sum_{i_1=1}^n x_{i_1}e_{i_1}\right).
		$$
		By bilinearity,
		$$
		v_1
		=
		\sum_{i_0,i_1=1}^n
		x_{i_0}x_{i_1} e_{i_0}e_{i_1}.
		$$
		Using the structure constants,
		$$
		e_{i_0}e_{i_1}
		=
		\sum_{k=1}^n \alpha_{i_0 i_1}^k e_k,
		$$
		we obtain
		$$
		v_1
		=
		\sum_{k=1}^n
		\left(
		\sum_{i_0,i_1=1}^n
		\alpha_{i_0 i_1}^k x_{i_0}x_{i_1}
		\right)e_k.
		$$
		Therefore
		$$
		v_1^k
		=
		\sum_{i_0,i_1=1}^n
		\alpha_{i_0 i_1}^k x_{i_0}x_{i_1}.
		$$
		
		In general, for $m\geq 1$, the coefficient of $e_k$ in $v_m$ is given by
		\begin{equation}\label{formulaexplicita}
			v_m^k
			=
			\sum_{i_0,\ldots,i_m=1}^n
			\sum_{a_1,\ldots,a_{m-1}=1}^n
			\alpha_{i_0 i_1}^{a_1}
			\alpha_{i_2 a_1}^{a_2}
			\alpha_{i_3 a_2}^{a_3}
			\cdots
			\alpha_{i_m a_{m-1}}^k
			x_{i_0}x_{i_1}\cdots x_{i_m}.
		\end{equation}
		For $m=1$, this formula is interpreted as
		$$
		v_1^k
		=
		\sum_{i_0,i_1=1}^n
		\alpha_{i_0 i_1}^k x_{i_0}x_{i_1},$$
		
		because there are no intermediate indices $a_1,\ldots,a_{m-1}$.
		
		Let us explain the meaning of the indices in the general formula. The indices
		$$
		i_0,i_1,\ldots,i_m
		$$
		come from the $m+1$ occurrences of $x$ in the right-nested product
		$$
		v_m=x(x(\cdots x(xx)\cdots)).
		$$
		The first multiplication is
		$$
		e_{i_0}e_{i_1},
		$$
		whose coefficient of $e_{a_1}$ is $\alpha_{i_0 i_1}^{a_1}$. The next multiplication is
		$$
		e_{i_2}e_{a_1},
		$$
		whose coefficient of $e_{a_2}$ is $\alpha_{i_2 a_1}^{a_2}$. Continuing in this way, after $m$ multiplications, the final basis vector is $e_k$. Thus the indices
		$$
		a_1,\ldots,a_{m-1}
		$$
		are intermediate indices, recording the basis vector obtained after each intermediate multiplication. It follows directly from the formula that each $v_m^k$ is homogeneous of degree $m+1$ in $x_1,\ldots,x_n$. Moreover, each coefficient is a polynomial in the structure constants. More precisely, each summand contains $m$ structure constants.\\
		
		Let $M=M(x)$ be the matrix with elements $M_{k,j}=v_{j-1}^k$, we have:
		$$
		\det M
		=
		\det
		\big(x,L_x(x),L_x^2(x),\ldots,L_x^{n-1}(x)\big).
		$$
		Using the Leibniz formula for the determinant, we get
		$$
		\det M=\sum_{\sigma\in S_n}
		\text{sgn}(\sigma)
		\prod_{m=1}^n M_{\sigma(m),m}
		=\sum_{\sigma\in S_n}
		\operatorname{sgn}(\sigma)
		\prod_{m=0}^{n-1}
		v_m^{\sigma(m+1)}.
		$$
		
		Since
		$$
		M_{\sigma(j),j}(x)=v_{j-1}^{\sigma(j)},
		$$
		we may rewrite this as
		$$
		\det M
		=
		\sum_{\sigma\in S_n}
		\operatorname{sgn}(\sigma)
		\prod_{m=0}^{n-1}
		v_m^{\sigma(m+1)}.
		$$
		For each $m$, the term $v_m^{\sigma(m+1)}$ is homogeneous of degree $m+1$ in $x_1,\ldots,x_n$. Therefore each summand in the determinant is homogeneous of degree
		$$
		1+2+\cdots+n=\frac{n(n+1)}{2}.
		$$
		Consequently, $\det M$ is a homogeneous polynomial of degree $N:=\frac{n(n+1)}{2}$ in the variables $x_1,\ldots,x_n$.\\
		
		Thus $\det M$ can be written in the form
		$$
		\det M
		=
		\sum_{\beta_1+\cdots+\beta_n=N}
		P_{(\beta_1,\ldots,\beta_n)}(\alpha)\,
		x_1^{\beta_1}\cdots x_n^{\beta_n}.
		$$
		For simplicity, we will denote $\beta:=(\beta_1,\ldots,\beta_n)$ and
		we also write $P_\beta(\alpha)=P_{(\beta_1,\ldots,\beta_n)}(\alpha)$. By Equation \ref{formulaexplicita}, each coefficient $P_\beta(\alpha)$ is a polynomial in the structure constants $\alpha_{ij}^k$.
		Since $K$ is algebraically closed, it is infinite. Hence $\det M(x)=0$ for every $x\in K^n$ if and only if every coefficient $P_\beta(\alpha)$ vanishes. The alicyclic laws therefore form a Zariski-closed subset. This subset is invariant under change of basis: an algebra isomorphism $h$ satisfies $h(x^j)=h(x)^j$, so it preserves the linear dependence of $x,\ldots,x^n$. If $A$ is alicyclic, this closed subset contains its orbit and consequently its orbit closure. Thus $A\to B$ implies that $B$ is alicyclic.
	\end{proof}
	From now on, let $A$ have dimension three. If $T$ is a linear operator on a three-dimensional space, we denote its characteristic polynomial by $\chi_T(\lambda)=\lambda^3+s_1\lambda^2+s_2\lambda+s_3$, $s_1, s_2, s_3 \in K$. In case $T=L_a$ for some element $a \in A$, we may write $s_i(a)$ instead of $s_i$ to emphasize this relation. We note that given an ordered basis $B$ for $A$ and $(a_1,a_2,a_3)=[a]_B = a\in A$ an arbitrary vector, the coefficients $s_i(a)$ are homogeneous polynomials in function of $a_1,a_2,a_3$. When the coordinates of $a$ are clear from the context, we may also write $s_i(a_1,a_2,a_3)$ to indicate that $s_i(a)$ also depends on the coordinates of $a\in A$.
	\begin{definition}
		Let $A$ be a three-dimensional algebra, we say that $A$ is $(\kappa,\alpha,\beta)$-characteristic, for $\alpha\neq\beta$,$\alpha,\beta\in\{1,2,3\},\kappa\in K$,if it satisfies:
		$$
		\forall a\in A\Bigg(  \chi_{L_a}(\lambda)=\lambda^3+s_1(a)\lambda^2+s_2(a)\lambda+s_3(a), s_\alpha(a)=\kappa s_\beta(a)\Bigg)
		$$
	\end{definition}
	
	That is, for every element $a$, its corresponding operator $L_a$ has a characteristic polynomial whose coefficients of its quadratic, linear or constant terms are a scalar multiple of some other coefficient, with the exception of the cubic term and this scalar is global, that is, it is the same for every $a\in A$. The following result is an immediate application of Proposition \ref{teohomopoly}.
	
	\begin{corollary}\label{principal}
		Let $A$ and $B$ be two algebras. If $A$ is $(\kappa,\alpha,\beta)$-characteristic and there is a deformation from $A$ to $B$, then $B$ is $(\kappa,\alpha,\beta)$-characteristic.
	\end{corollary}
	
	The relation $s_\alpha=\kappa s_\beta$ is a particular case of $s_\alpha=p(s_\beta)$ where $p$ is the homogeneous polynomial $p(y)=\kappa y$. This suggests the following generalization of which Corollary \ref{characteristic} follows by Proposition \ref{teohomopoly}
	\begin{definition}\label{defcharacteristic}
		We say that a three-dimensional algebra $A$ is $(\alpha,\beta)_p$-characteristic  for $\alpha\neq\beta$,$\alpha,\beta\in\{1,2,3\}$ and $p(x)=\kappa x^n$, where $n\geq1$ and $\kappa\in K$ if:
		$$
		\forall a\in A\Bigg( \chi_{L_a}(\lambda)=\lambda^3+s_1(a)\lambda^2+s_2(a)\lambda+s_3(a), s_\alpha(a)=p(s_\beta(a))\Bigg).
		$$
	\end{definition}
	\begin{corollary}\label{characteristic}
		Let $A$ and $B$ be two algebras. If $A$ is $(\alpha,\beta)_p$-characteristic and there is a deformation from $A$ to $B$, then $B$ is $(\alpha,\beta)_p$-characteristic.
	\end{corollary}
	For both corollaries, expand the stated identity as a polynomial in the coordinates of $a$. Its coefficients are polynomials in the structure constants by Proposition \ref{teohomopoly}. Their simultaneous vanishing defines a Zariski-closed subset. It is invariant under isomorphisms because multiplication operators are conjugate after transport of the element $a$. Consequently the identity persists in orbit closures. The scalar $\kappa$ and the polynomial $p$ are fixed throughout each assertion. In particular, different homogeneous degrees in Definition \ref{defcharacteristic} may force both coefficients to vanish; this does not affect the preservation argument.
	
	We present a proposition, concerning polynomial identities and degenerations of algebras that will also be used in our classification.
	
	\begin{proposition}
		Let $K$ be an arbitrary field and $V$ be an $n$-dimensional $K$-vector space. Let $P(x_1,\ldots,x_k)$ be a nonassociative homogeneous polynomial of degree at least two. The set
		$$
		V(P)=\{A\in K^{n^3}: P_A(a_1,\ldots,a_k)=0,\forall a_i\in A\},
		$$
		is a Zariski-closed set of $K^{n^3}$.
	\end{proposition}
	\begin{proof}
		Fix a basis $\{e_1,\ldots,e_n\}$ for $V$ and identify an algebra $A$ with its structure constants $(\alpha_{i,j}^k)$, we abuse notation and denote this set by $A$. Fix an element $(a_1,\ldots,a_k)\in V^k$. Consider the map
		$$
		\widetilde{P}_a\colon K^{n^3}\to V,\qquad A\mapsto P_A(a_1,\ldots,a_k).
		$$
		
		If $P$ has degree $d$, each coordinate of $\widetilde{P}_a$ is a homogeneous polynomial of degree $d-1$ in the structure constants of $A$: every degree-$d$ monomial uses $d-1$ multiplications. Consequently, the set
		$$
		Z_a=\{ A\in K^{n^3}: \widetilde{P}_a(A)=0 \},
		$$
		is Zariski-closed. Since $V(P)=\displaystyle\bigcap_{a\in V^k} Z_a$ is the intersection of closed sets, it is also closed, proving the statement.
	\end{proof}
	The following criteria will be used for the determination of degenerations. All comparisons of automorphism dimensions concern nonisomorphic algebras and use the dimensions in the tables.
	\begin{itemize}
		\item $\underset{Aut}{\nrightarrow}$: If $A$ properly degenerates to $B$, then $\dim\operatorname{Aut}(A)<\dim\operatorname{Aut}(B)$, since the boundary of an algebraic group orbit consists of orbits of strictly smaller dimension [\cite{aut}, Section 8.3].
		\item $\underset{Ass}{\nrightarrow}$: Associativity is a polynomial identity, so an associative algebra cannot degenerate to a nonassociative algebra.
		\item $\underset{Unital}{\nrightarrow}$: In the entries carrying this label, $\det L_x=0$ for every $x$ in the source, whereas the target has an identity $1$ and $\det L_1=1$. The determinant condition is polynomial and invariant under isomorphism, and hence persists under degeneration.
		\item $\underset{Alicyclic}{\nrightarrow}$: The source is alicyclic and the target is not; apply Proposition \ref{propositionAlicylic}.
		\item Labels $x^2=0$, $x^2y=0$ and $x^4=0$ indicate polynomial identities satisfied by the source and violated by the target.
		\item $\underset{L_{x^2}=c(x)I}{\nrightarrow}$: Every $L_{x^2}$ in the source is scalar, but this fails in the target. Vanishing of the off-diagonal entries and of differences of diagonal entries gives polynomial equations for this condition.
		\item $\underset{trL_x}{\nrightarrow}$: The source satisfies $\operatorname{tr}L_x=0$ for every $x$, whereas the target does not.
		\item $\underset{tr(L_{xy})}{\nrightarrow}$: The rank of the bilinear form $B_A(x,y)=\operatorname{tr}(L_{xy})$ cannot increase in a degeneration, since the rank bounds are given by minors of the matrix that defines the form.
		\item $\underset{Ann}{\nrightarrow}$: $\dim\operatorname{Ann}(A)>\dim\operatorname{Ann}(B)$. The annihilator is the kernel of a matrix polynomial in the structure constants, so its dimension cannot decrease under degeneration.
		\item $\underset{Structure}{\nrightarrow}$: For $J_{3,6}$ write
		$x=s(x)1+r(x)$ with $r(x)$ in its square-zero radical. Then
		$xy=s(x)y+s(y)x+s(x)s(y)1$. In an arbitrary ordered basis
		$\{e_1,e_2,e_3\}$, put $s_i=s(e_i)$ and write $1=\sum_i u_i e_i$.
		The structure constants satisfy
		$$
		c_{11}^3=s_1^2u_3,\qquad c_{22}^3=s_2^2u_3,
		\qquad c_{12}^3=s_1s_2u_3.
		$$
		Thus $c_{11}^3c_{22}^3=(c_{12}^3)^2$ on the entire orbit and hence
		on its closure. In the displayed table of $J_{3,25}$, these three
		constants are $0,0,1$, respectively, so the equation fails.
		Consequently $J_{3,6}\nrightarrow J_{3,25}$.
		\item $\underset{dimZ}{\nrightarrow}$: $\dim Z(A)>\dim Z(B)$, where $Z(A)$ is the associative center: its elements commute with all elements and associate in all three positions. In the present commutative conditions,
		$$
		Z(A)=\{z:(z,x,y)=(x,z,y)=(x,y,z)=0\text{ for all }x,y\in A\},
		$$
		where $(x,y,z)=(xy)z-x(yz)$. These are linear equations in $z$, with polynomial coefficients in the structure constants; hence $\dim Z$ cannot decrease under degeneration (we refer to \cite[Proposition 4.58]{MM}, the proof remains unchanged if ${\rm char}(K)=2)$).
		\item $\underset{Characteristic}{\nrightarrow}$:Use a fixed polynomial identity among the coefficients $s_i(x)$ of the characteristic polynomial of the operator $L_x$, as in Corollary \ref{characteristic}. The specific identities used below are $s_2=s_1^2$ and $s_1^2=\gamma^2s_2$.
		\item $\underset{Det}{\nrightarrow}$:
		For a three-dimensional commutative algebra over a field of characteristic two satisfying
		$x^2=0$, with basis $\{e_1,e_2,e_3\}$. Let $T_A=[e_2e_3\mid e_3e_1\mid e_1e_2]$ be the square matrix whose entries are the coordinates of these vectors.
		For any fixed $\lambda\in K$, the condition
		$$
		\det(T_A+\lambda T_A^{\mathsf T})=0
		$$
		is preserved under degeneration, since the determinant
		is polynomial in the structure constants and transforms
		by multiplication by $(\det g)^{-1}$ under a change
		of basis. Taking $\lambda$ such that $1+\lambda^2+\gamma^2\lambda=0$ proves that there is no deformation $J_{3,31}^{\gamma \neq0}\to J_{3,31}^{0}$.
		\item  $\underset{dim A^2}{\nrightarrow}$, there is no degeneration $A\to B$ and this follows by the relation $\dim A^2<\dim B^2$.
		\item $\underset{trivial}{\rightarrow}$: The basis $f_i=te_i$ gives degeneration to the zero algebra.
		\item $\underset{t}{\rightarrow}$: The degeneration follows by transitivity.
		\item ${\rightarrow}$: The source and target coincide, or an explicit basis family is listed below. Parameters remain fixed as $t$ varies.
	\end{itemize}

	\subsection{Two-dimensional algebras}
	
	\begin{theorem}\label{classification2geometry}
		For the variety $\rm{Jor}_2$ over an algebraically closed field of characteristic two, we have:
		\begin{itemize}
			\item The algebras $J_{2,5}^0$ and $J_{2,7}$ are rigid;
			\item The variety has two irreducible components given by $\overline{(J_{2,5}^0)^ {GL(V)}}$ and $\overline{J_{2,7}^ {GL(V)}}$;
			\item Write $\mathcal O(A)=A^{GL(V)}$. The degeneration relations yield the following orbit closures:
			$$
			\begin{aligned}
				\overline{\mathcal O(J_{2,1})}
				&=\mathcal O(J_{2,1}),\\[2mm]
				\overline{\mathcal O(J_{2,2})}
				&=\mathcal O(J_{2,2})\cup\mathcal O(J_{2,1}),\\[2mm]
				\overline{\mathcal O(J_{2,3})}
				&=\mathcal O(J_{2,3})\cup\mathcal O(J_{2,2})
				\cup\mathcal O(J_{2,1}),\\[2mm]
				\overline{\mathcal O(J^0_{2,4})}
				&=\mathcal O(J^0_{2,4})\cup\mathcal O(J_{2,2})
				\cup\mathcal O(J_{2,1}),\\[2mm]
				\overline{\mathcal O(J^0_{2,5})}
				&=\mathcal O(J^0_{2,5})\cup\mathcal O(J^0_{2,4})\cup\mathcal O(J_{2,3}) 
				\cup\mathcal O(J_{2,2})\cup\mathcal O(J_{2,1}),\\[2mm]
				\overline{\mathcal O(J_{2,7})}
				&=\mathcal O(J_{2,7})\cup\mathcal O(J_{2,1}).
			\end{aligned}
			$$
		\end{itemize}
	\end{theorem}
	
	The statement of the above theorem will be proven by an explicit computation of the isomorphism classes and degenerations, when they exist. The isomorphism classes, their change of bases and the dimension of the automorphism groups are the subject of the following table.
	\begin{table}[!h]
		\begin{adjustbox}{max totalsize={\textwidth}{0.88\textheight},center}
			\begin{tabular}{|l|l|l|l|l|}
				\hline
				Algebra $J$ & Change of basis & Table & ${\rm dim}({\rm Aut}(J))$ & Notes \\ \hline
				$J_{2,1}$    &      -           &   $\langle a,b \ |\rangle$    &     4   & Zero-algebra      \\ \hline
				$J_{2,2}$    &        -         &    $\langle a,a^2\ | \rangle$    &     2   &   Associative   \\ \hline
				$J_{2,3}$    &           -      &   $\langle e,a\ |\ e^2=e \rangle$    &   1     &  Associative, non-unital     \\ \hline
				$J^\rho_{2,4}$    &      $\{1,a+(\sqrt{\rho})1\}$           &   $J^0_{2,4}=\langle 1,a\ |\rangle$    &    1    &     Associative, unital  \\ \hline
				$J^\kappa_{2,5}$    &   $\{1,a+\beta 1\};\beta^2+\kappa=\beta$              &    $J^0_{2,5}=\langle 1, a\ |\ a^2=a\rangle $   &   0     &     Associative, unital  \\ \hline
				$J_{2,7}$    &   -              &    $J_{2,7}=\langle a, b\ |\ ab=a \rangle $   &   2     &    Non-associative Jordan  \\ \hline
			\end{tabular}
		\end{adjustbox}
		\caption{Two-dimensional algebras over an algebraically closed field of characteristic 2.}
	\end{table}
	
	From where we extract the table describing their degenerations.
	\begin{center}
		\begin{tabular}{|c|c|c|c|c|c|c|}
			\hline 
			Algebra & $J_{2,1}$ & $J_{2,2}$ & $J_{2,3}$ & $J^{0}_{2,4}$ & $J^0_{2,5}$ & $J_{2,7}$\tabularnewline
			\hline 
			\hline 
			$J_{2,1}$ & $\shortrightarrow$ & $\underset{Aut}{\nrightarrow}$ & $\underset{Aut}{\nrightarrow}$ & $\underset{Aut}{\nrightarrow}$ & $\underset{Aut}{\nrightarrow}$ & $\underset{Aut}{\nrightarrow}$\tabularnewline
			\hline 
			$J_{2,2}$ & $\shortrightarrow$ & $\shortrightarrow$ & $\underset{Aut}{\nrightarrow}$ & $\underset{Aut}{\nrightarrow}$ & $\underset{Aut}{\nrightarrow}$ & $\underset{Aut}{\nrightarrow}$\tabularnewline
			\hline 
			$J_{2,3}$ & $\overset{\longrightarrow}{t}$ & $\shortrightarrow$ & $\shortrightarrow$ & $\underset{Aut}{\nrightarrow}$ & $\underset{Aut}{\nrightarrow}$ & $\underset{Ass}{\nrightarrow}$\tabularnewline
			\hline 
			$J^0_{2,4}$ & $\overset{\longrightarrow}{t}$ & $\shortrightarrow$ & $\underset{Aut}{\nrightarrow}$ & $\shortrightarrow$ & $\underset{Aut}{\nrightarrow}$ & $\underset{Ass}{\nrightarrow}$\tabularnewline
			\hline 
			$J^0_{2,5}$ & $\overset{\longrightarrow}{t}$ & $\overset{\longrightarrow}{t}$ & $\shortrightarrow$ & $\shortrightarrow$ & $\shortrightarrow$ & $\underset{Ass}{\nrightarrow}$\tabularnewline
			\hline 
			$J_{2,7}$ & $\shortrightarrow$ & $\underset{Aut}{\nrightarrow}$ & $\underset{Aut}{\nrightarrow}$ & $\underset{Aut}{\nrightarrow}$ & $\underset{Aut}{\nrightarrow}$ & $\shortrightarrow$\tabularnewline
			\hline 
		\end{tabular}
	\end{center}
	
	\begin{table}[h!]
		\begin{adjustbox}{max totalsize={\textwidth}{0.88\textheight},center}
			\begin{tabular}{|l|l|}
				\hline
				Deformation & Change of basis \\ \hline
				$J_{2,3}\to J_{2,2}$   &   $f_1=te+a,\ f_2=t^2e$              \\ \hline
				$J^0_{2,4}\to J_{2,2}$   &     $f_1=t+a,\ f_2=t^2$            \\ \hline
				$J^0_{2,5}\to J_{2,3}$ & $f_1=1+a,\ f_2=ta$ \\ \hline
				$J^0_{2,5} \to J_{2,4}^0$   &      $f_1=1,\ f_2=ta$
				\\ \hline
			\end{tabular}
		\end{adjustbox}
	\end{table}
	\begin{proof}[Proof of Theorem \ref{classification2geometry}]
		Every orbit belongs to the closure of $J_{2,5}^0$ or $J_{2,7}$, and neither of these two orbits lies in the closure of the other. Their closures are therefore the two irreducible components. Each of these orbits is open: its complement is the finite union of the closures of the other orbits. Every other orbit has a proper incoming degeneration and hence is not open.
	\end{proof}
	\subsection{Three-dimensional algebras}
	\begin{theorem}\label{zeropotentcomponent}
		Let $\dim V=3$ and let $K$ be algebraically closed of characteristic two. The subset
		$$
		\mathcal Z=\{A\in\mathrm{Jor}_3:x^2=0\text{ for every }x\in A\}
		$$
		is closed and irreducible, and is isomorphic to affine space $K^9$.
	\end{theorem}
	\begin{proof}
		For a commutative bilinear product in characteristic two, $x^2=0$ for every $x$ is equivalent to $e_i^2=0$ for the three basis vectors. The three mixed products $e_1e_2,e_1e_3,e_2e_3$ are arbitrary vectors in $V$, giving nine independent scalar coordinates. Every such product satisfies the Jordan identity because $x^2=0$. Thus $\mathcal Z$ is a linear closed subspace isomorphic to $K^9$, and is irreducible.
	\end{proof}
	\clearpage

	\begin{table}[!h]
		\begin{adjustbox}{max totalsize={\textwidth}{0.88\textheight},center}
			\begin{tabular}{|l|l|l|l|}
				\hline
				$Algebra$ & $Change\ of\ basis$ &  $Class$ & dim$Aut(J)$ \\ \hline
				$J_{3,1}$ & $-$ & $J_{3,1}=\lan a,b,c\mid\ran$ &  9 \\ \hline
				$J_{3,2}$ & $-$ & $J_{3,2}=\lan a,b,c\mid a^ 2=b\ran$ &  5  \\ \hline
				$J_{3,3}^{\delta\neq 0} $ & $\{a,b/\sqrt{\delta},c\}$ & $J_{3,3}^1=\lan a,b,c\mid a^2=c,b^2= c\ran$ & 4   \\ \hline
				$J_{3,3}^0$& $\{a,c,b\}$ &  $J_{3,2}$& -\\ \hline
				$J_{3,4}$& $-$ & $J_{3,4}=\lan a,b,c\mid a^2=b,ab=c \ran$ & 3     \\ \hline
				$J_{3,5}$& $-$ & $J_{3,5}=\lan 1,a,b\mid b^2=b\ran$   & 1   \\ \hline
				$J_{3,6}$& $-$ &$J_{3,6}=\lan 1,a,b\ran$  & 4  \\ \hline
				$J_{3,7}$& $-$ & $J_{3,7}=\lan 1,a,a^2\ran$   & 2   \\ \hline
				$J_{3,8}$& $-$ &$ J_{3,8}=\lan a,b,c\mid a^2=b,c^2=c\ran$  & 2    \\ \hline
				$J_{3,9}$& $-$ & $J_{3,9}=\lan a,b,c\mid a^ 2=a \ran$   & 4  \\ \hline
				$J_{3,10}^{\rho\neq 0}$& $\{a,a+b/\sqrt{\rho},c\}$ & $J_{3,10}^0=\lan a,b,c\mid a^2=a,ab=b\ran$   & - \\ \hline
				$J_{3,10}^0$& $-$ & $J_{3,10}^0$ & 2  \\ \hline
				$J_{3,11}^{\kappa \neq0}$& $\{a,b+ua,c\},u^2=u+\kappa$ & $J_{3,11}^0=\lan a,b,c\mid a^2=a,ab=b,b^2=b \ran$  & -\\ \hline
				$J_{3,11}^0$& $-$ &  $J_{3,11}^0$& 1  \\ \hline
				$J_{3,12}^{\rho\neq 0}$& $\{a,b,c+\sqrt{\rho}b\}$ &$J_{3,12}^0=\langle a,b,c\mid a^2=a,b^2=b,bc=c\rangle$   & -\\ \hline
				$J_{3,12}^0$& $\{a+b,c,a\}$ & $J_{3,5}$  & -\\ \hline
				$J_{3,13}^{\kappa \neq 0}$& $\{a,b+ua,c\},\ u^2+u=\kappa$ & $J_{3,13}^0=\langle a,b,c\mid a^2=a,b^2=b,c^2=c,ab=b\rangle$  & -\\ \hline
				$J_{3,13}^0$& $-$ & $J_{3,13}^0$   & 0  \\ \hline
				$J_{3,15}^{\ell\neq 0}$& $\{a,a^2,b/\sqrt{\ell}\}$ & $J_{3,15}^1=\lan a,a^2,b\mid ab=b, b^2=a^2\ran$ & 2 \\ \hline
				$J_{3,15}^0$& $-$ & $J_{3,15}^0=\lan a,a^2,b\mid ab=b\ran$   & 3  \\ \hline
				$J_{3,16}^{\alpha\neq 0,0}$& $\{1,a/\sqrt{\alpha},b\}$ &$J_{3,16}^{1,0}=\lan 1,a,b \mid a^2=1\ran$ & 2   \\ \hline
				$J_{3,16}^{0,0}$& $-$ &$J_{3,6}$ & - \\ \hline		
				$J_{3,16}^{0,\beta\neq 0}$& $\{1,b/\sqrt{\beta},a\}$ & $J_{3,16}^{1,0}$ & - \\ \hline			
				$J_{3,16}^{\alpha\neq0,\beta\neq0}$& $\{1,a/\sqrt{\alpha},b/\sqrt{\beta}\}$ &$J_{3,16}^{1,1}=\lan 1,a,b\mid a^2=1,\ b^2=1\ran$  & 1  \\ \hline		
				$J_{3,17}$& $-$ &$J_{3,17}=\lan 1,a,b\mid ab=1\ran$ & 3  \\ \hline
				$J_{3,24}$& $-$ & $J_{3,24}=\langle a,b,c\mid a^2=a,bc=b\rangle$  & 2  \\ \hline
				$J_{3,25}$& $-$ & $J_{3,25}=\lan a,b,c\mid ab=c \ran$ &  6    \\ \hline
				$J_{3,26}$& $-$ &$J_{3,26}=\lan a,b,c \mid ab=c,bc=c\ran$ & 4   \\ \hline
				$J_{3,27}$&  $-$ &$J_{3,27}=\lan a,b,c \mid ac=a,bc=b\ran$ & 6   \\ \hline
				$J_{3,28}^{\gamma\neq 0}$&  $-$ &$J_{3,28}^{\gamma\neq 0}=\lan a,b,c \mid ab=c,ac=b+\gamma c\ran$ & 4   \\ \hline
				$J_{3,28}^0$&  $-$ &$J_{3,28}^0=\lan a,b,c \mid ab=c,ac=b\ran$ & 4   \\ \hline
				$J_{3,29}$&  $-$ &$J_{3,29}=\lan a,b,c \mid ab=c,bc=b\ran$ & 1    \\ \hline
				$J_{3,30}$&  $-$ &$J_{3,30}=\lan a,b,c \mid ab=c,ac=c,bc=b+c\ran$ & 2  \\ \hline
				$J_{3,31}^{\gamma\neq 0}$&  $-$ &$J_{3,31}^{\gamma\neq 0}=\lan a,b,c \mid ab=c,ac=b,bc=a+\gamma b\ran$ & 1   \\
				\hline
				$J_{3,31}^{0}$&  $-$ &$J_{3,31}^0=\lan a,b,c \mid ab=c,ac=b,bc=a\ran$ & 3  \\\hline
				$J_{3,32}$&  $-$ &$J_{3,32}=\lan a,b,c \mid ab=c,ac=b,bc=a+b+c\ran$ &  1  \\ \hline
			\end{tabular}
		\end{adjustbox}
		\caption{Isomorphism classes over an algebraically closed field of characteristic two}
	\end{table}

	In the following matrix, the parameterized rows and columns refer to
	$\gamma\neq0$. A diagonal arrow refers to the same isomorphism class;
	it does not assert a degeneration between distinct members of a family.
	Parameters in different families are independent. For the basis listed
	for $J_{3,31}^{\gamma}\to J_{3,28}^{\delta}$, set
	$\beta=\delta/\gamma$ and $r=1+\gamma/\delta$, where $\gamma\delta\neq0$.
	Its determinant is $\beta^2t^2$, and its nonzero products are
	$$
	f_1f_2=f_3,\qquad
	f_1f_3=f_2+\delta f_3+\beta rt f_1,\qquad
	f_2f_3=t^2f_1+\beta rt f_2.
	$$
	These specialize to $J_{3,28}^{\delta}$.

	\clearpage
	\begin{table}
		\adjustbox{max totalsize={\textwidth}{1\textheight}}{
			\begin{adjustbox}{angle=90}
				\begin{tabular}{|c|c|c|c|c|c|c|c|c|c|c|c|c|c|c|c|c|c|c|c|c|c|c|c|c|c|c|c|c|}
					\hline 
					Algebra & $J_{3,1}$ & $J_{3,2}$ & $J_{3,3}^{1}$ & $J_{3,4}$ & $J_{3,5}$ & $J_{3,6}$ & $J_{3,7}$ & $J_{3,8}$ & $J_{3,9}$ & $J_{3,10}^{0}$ & $J_{3,11}^{0}$ & $J_{3,13}^{0}$ & $J_{3,15}^1$ & $J_{3,15}^{0}$ & $J_{3,16}^{1,0}$ & $J_{3,16}^{1,1}$ & $J_{3,17}$ & $J_{3,24}$ & $J_{3,25}$ & $J_{3,26}$ & $J_{3,27}$ & $J_{3,28}^{\gamma\neq0}$ & $J_{3,29}$ & $J_{3,30}$ & $J_{3,31}^{\gamma\neq0}$ & $J_{3,32}$ & $J_{3,28}^{0}$ & $J_{3,31}^{0}$ \tabularnewline
					\hline 
					\hline 
					$J_{3,1}$ & $\longrightarrow$ & $\underset{Aut}{\nrightarrow}$ & $\underset{Aut}{\nrightarrow}$ & $\underset{Aut}{\nrightarrow}$ & $\underset{Aut}{\nrightarrow}$ & $\underset{Aut}{\nrightarrow}$ & $\underset{Aut}{\nrightarrow}$ & $\underset{Aut}{\nrightarrow}$ & $\underset{Aut}{\nrightarrow}$ & $\underset{Aut}{\nrightarrow}$ & $\underset{Aut}{\nrightarrow}$ & $\underset{Aut}{\nrightarrow}$ & $\underset{Aut}{\nrightarrow}$ & $\underset{Aut}{\nrightarrow}$ & $\underset{Aut}{\nrightarrow}$ & $\underset{Aut}{\nrightarrow}$ & $\underset{Aut}{\nrightarrow}$ & $\underset{Aut}{\nrightarrow}$ & $\underset{Aut}{\nrightarrow}$ & $\underset{Aut}{\nrightarrow}$ & $\underset{Aut}{\nrightarrow}$ & $\underset{Aut}{\nrightarrow}$ & $\underset{Aut}{\nrightarrow}$ & $\underset{Aut}{\nrightarrow}$ & $\underset{Aut}{\nrightarrow}$ & $\underset{Aut}{\nrightarrow}$ &  $\underset{Ass}{\nrightarrow}$ &  $\underset{Ass}{\nrightarrow}$\tabularnewline
					\hline 
					$J_{3,2}$ & $\underset{trivial}{\longrightarrow}$ & $\longrightarrow$ & $\underset{Aut}{\nrightarrow}$ & $\underset{Aut}{\nrightarrow}$ & $\underset{Aut}{\nrightarrow}$ & $\underset{Aut}{\nrightarrow}$ & $\underset{Aut}{\nrightarrow}$ & $\underset{Aut}{\nrightarrow}$ & $\underset{Aut}{\nrightarrow}$ & $\underset{Aut}{\nrightarrow}$ & $\underset{Aut}{\nrightarrow}$ & $\underset{Aut}{\nrightarrow}$ & $\underset{Aut}{\nrightarrow}$ & $\underset{Aut}{\nrightarrow}$ & $\underset{Aut}{\nrightarrow}$ & $\underset{Aut}{\nrightarrow}$ & $\underset{Aut}{\nrightarrow}$ & $\underset{Aut}{\nrightarrow}$ & $\underset{Ann}{\nrightarrow}$ & $\underset{Aut}{\nrightarrow}$ & $\underset{Ass}{\nrightarrow}$ & $\underset{Aut}{\nrightarrow}$ & $\underset{Aut}{\nrightarrow}$ & $\underset{Aut}{\nrightarrow}$ & $\underset{Aut}{\nrightarrow}$ & $\underset{Aut}{\nrightarrow}$ & $\underset{Ass}{\nrightarrow}$ &  $\underset{Ass}{\nrightarrow}$\tabularnewline
					\hline 
					$J_{3,3}^{1}$ & $\underset{trivial}{\longrightarrow}$ & $\longrightarrow$ & $\longrightarrow$ & $\underset{Aut}{\nrightarrow}$ & $\underset{Aut}{\nrightarrow}$ & $\underset{Aut}{\nrightarrow}$ & $\underset{Aut}{\nrightarrow}$ & $\underset{Aut}{\nrightarrow}$ & $\underset{Aut}{\nrightarrow}$ & $\underset{Aut}{\nrightarrow}$ & $\underset{Aut}{\nrightarrow}$ & $\underset{Aut}{\nrightarrow}$ & $\underset{Aut}{\nrightarrow}$ & $\underset{Aut}{\nrightarrow}$ & $\underset{Aut}{\nrightarrow}$ & $\underset{Aut}{\nrightarrow}$ & $\underset{Aut}{\nrightarrow}$ & $\underset{Aut}{\nrightarrow}$ & $\longrightarrow$ & $\underset{Aut}{\nrightarrow}$ & $\underset{Ass}{\nrightarrow}$ & $\underset{Aut}{\nrightarrow}$ & $\underset{Aut}{\nrightarrow}$ & $\underset{Aut}{\nrightarrow}$ & $\underset{Aut}{\nrightarrow}$ & $\underset{Aut}{\nrightarrow}$  & $\underset{Ass}{\nrightarrow}$ &  $\underset{Ass}{\nrightarrow}$\tabularnewline
					\hline 
					$J_{3,4}$ & $\underset{trivial}{\longrightarrow}$ & $\longrightarrow$ & $\longrightarrow$ & $\longrightarrow$ & $\underset{Aut}{\nrightarrow}$ & $\underset{Unital}{\nrightarrow}$ & $\underset{Aut}{\nrightarrow}$ & $\underset{Aut}{\nrightarrow}$ & $\underset{x^4=0}{\nrightarrow}$ & $\underset{Aut}{\nrightarrow}$ & $\underset{Aut}{\nrightarrow}$ & $\underset{Aut}{\nrightarrow}$ & $\underset{Aut}{\nrightarrow}$ & $\underset{Aut}{\nrightarrow}$ & $\underset{Unital}{\nrightarrow}$ & $\underset{Aut}{\nrightarrow}$ & $\underset{Aut}{\nrightarrow}$ & $\underset{Aut}{\nrightarrow}$ & $\longrightarrow$ & $\underset{Ass}{\nrightarrow}$ & $\underset{Ass}{\nrightarrow}$ & $\underset{Ass}{\nrightarrow}$ & $\underset{Aut}{\nrightarrow}$ & $\underset{Aut}{\nrightarrow}$ & $\underset{Aut}{\nrightarrow}$ & $\underset{Aut}{\nrightarrow}$  & $\underset{Ass}{\nrightarrow}$ &  $\underset{Ass}{\nrightarrow}$\tabularnewline
					\hline 
					$J_{3,5}$ & $\underset{trivial}{\longrightarrow}$ & $\underset{t}{\longrightarrow}$ & $\underset{t}{\longrightarrow}$ & $\underset{t}{\longrightarrow}$ & $\longrightarrow$ & $\underset{t}{\longrightarrow}$ & $\longrightarrow$ & $\longrightarrow$ & $\underset{t}{\longrightarrow}$ & $\longrightarrow$ & $\underset{Aut}{\nrightarrow}$ & $\underset{Aut}{\nrightarrow}$ & $\underset{Ass}{\nrightarrow}$ & $\underset{Ass}{\nrightarrow}$ & $\underset{Ass}{\nrightarrow}$ & $\underset{Aut}{\nrightarrow}$ & $\underset{Ass}{\nrightarrow}$ & $\underset{Ass}{\nrightarrow}$ & $\underset{t}{\longrightarrow}$ & $\underset{Ass}{\nrightarrow}$ & $\underset{Ass}{\nrightarrow}$ & $\underset{Ass}{\nrightarrow}$ & $\underset{Aut}{\nrightarrow}$ & $\underset{Ass}{\nrightarrow}$ & $\underset{Aut}{\nrightarrow}$ & $\underset{Aut}{\nrightarrow}$  & $\underset{Ass}{\nrightarrow}$ &  $\underset{Ass}{\nrightarrow}$\tabularnewline
					\hline 
					$J_{3,6}$ & $\underset{trivial}{\longrightarrow}$ & $\longrightarrow$ & $\underset{Aut}{\nrightarrow}$ & $\underset{Aut}{\nrightarrow}$ & $\underset{Aut}{\nrightarrow}$ & $\longrightarrow$ & $\underset{Aut}{\nrightarrow}$ & $\underset{Aut}{\nrightarrow}$ & $\underset{Aut}{\nrightarrow}$ & $\underset{Aut}{\nrightarrow}$ & $\underset{Aut}{\nrightarrow}$ & $\underset{Aut}{\nrightarrow}$ & $\underset{Aut}{\nrightarrow}$ & $\underset{Aut}{\nrightarrow}$ & $\underset{Aut}{\nrightarrow}$ & $\underset{Aut}{\nrightarrow}$ & $\underset{Aut}{\nrightarrow}$ & $\underset{Aut}{\nrightarrow}$ & $\underset{Structure}{\nrightarrow}$ & $\underset{Aut}{\nrightarrow}$ & $\underset{Ass}{\nrightarrow}$ & $\underset{Aut}{\nrightarrow}$ & $\underset{Aut}{\nrightarrow}$ & $\underset{Aut}{\nrightarrow}$ & $\underset{Aut}{\nrightarrow}$ & $\underset{Aut}{\nrightarrow}$  & $\underset{Ass}{\nrightarrow}$ &  $\underset{Ass}{\nrightarrow}$\tabularnewline
					\hline 
					$J_{3,7}$ & $\underset{trivial}{\longrightarrow}$ & $\longrightarrow$ & $\longrightarrow$ & $\longrightarrow$ & $\underset{Aut}{\nrightarrow}$ & $\longrightarrow$ & $\longrightarrow$ & $\underset{Aut}{\nrightarrow}$ & $\underset{Characteristic}{\nrightarrow}$ & $\underset{Aut}{\nrightarrow}$ & $\underset{Aut}{\nrightarrow}$ & $\underset{Aut}{\nrightarrow}$ & $\underset{Aut}{\nrightarrow}$ & $\underset{Ass}{\nrightarrow}$ & $\underset{Ass}{\nrightarrow}$ & $\underset{Aut}{\nrightarrow}$ & $\underset{Ass}{\nrightarrow}$ & $\underset{Aut}{\nrightarrow}$ & $\underset{t}{\longrightarrow}$ & $\underset{Ass}{\nrightarrow}$ & $\underset{Ass}{\nrightarrow}$ & $\underset{Ass}{\nrightarrow}$ & $\underset{Aut}{\nrightarrow}$ & $\underset{Aut}{\nrightarrow}$ & $\underset{Aut}{\nrightarrow}$ & $\underset{Aut}{\nrightarrow}$  & $\underset{Ass}{\nrightarrow}$ &  $\underset{Ass}{\nrightarrow}$\tabularnewline
					\hline 
					$J_{3,8}$ & $\underset{trivial}{\longrightarrow}$ & $\underset{t}{\longrightarrow}$ & $\underset{t}{\longrightarrow}$ & $\longrightarrow$ & $\underset{Aut}{\nrightarrow}$ & $\underset{Unital}{\nrightarrow}$ & $\underset{Aut}{\nrightarrow}$ & $\longrightarrow$ & $\longrightarrow$ & $\underset{Aut}{\nrightarrow}$ & $\underset{Aut}{\nrightarrow}$ & $\underset{Aut}{\nrightarrow}$ & $\underset{Aut}{\nrightarrow}$ & $\underset{Ass}{\nrightarrow}$ & $\underset{Ass}{\nrightarrow}$ & $\underset{Aut}{\nrightarrow}$ & $\underset{Ass}{\nrightarrow}$ & $\underset{Aut}{\nrightarrow}$ & $\underset{t}{\longrightarrow}$ & $\underset{Ass}{\nrightarrow}$ & $\underset{Ass}{\nrightarrow}$ & $\underset{Ass}{\nrightarrow}$ & $\underset{Aut}{\nrightarrow}$ & $\underset{Aut}{\nrightarrow}$ & $\underset{Aut}{\nrightarrow}$ & $\underset{Aut}{\nrightarrow}$  & $\underset{Ass}{\nrightarrow}$ &  $\underset{Ass}{\nrightarrow}$\tabularnewline
					\hline 
					$J_{3,9}$ & $\underset{trivial}{\longrightarrow}$ & $\longrightarrow$ & $\underset{Aut}{\nrightarrow}$ & $\underset{Aut}{\nrightarrow}$ & $\underset{Aut}{\nrightarrow}$ & $\underset{Aut}{\nrightarrow}$ & $\underset{Aut}{\nrightarrow}$ & $\underset{Aut}{\nrightarrow}$ & $\longrightarrow$ & $\underset{Aut}{\nrightarrow}$ & $\underset{Aut}{\nrightarrow}$ & $\underset{Aut}{\nrightarrow}$ & $\underset{Aut}{\nrightarrow}$ & $\underset{Aut}{\nrightarrow}$ & $\underset{Aut}{\nrightarrow}$ & $\underset{Aut}{\nrightarrow}$ & $\underset{Aut}{\nrightarrow}$ & $\underset{Aut}{\nrightarrow}$ & $\underset{Ann}{\nrightarrow}$ & $\underset{Aut}{\nrightarrow}$ & $\underset{Ass}{\nrightarrow}$ & $\underset{Aut}{\nrightarrow}$ & $\underset{Aut}{\nrightarrow}$ & $\underset{Aut}{\nrightarrow}$ & $\underset{Aut}{\nrightarrow}$ & $\underset{Aut}{\nrightarrow}$  & $\underset{Ass}{\nrightarrow}$ &  $\underset{Ass}{\nrightarrow}$\tabularnewline
					\hline 
					$J_{3,10}^{0}$ & $\underset{trivial}{\longrightarrow}$ & $\underset{t}{\longrightarrow}$ & $\underset{t}{\longrightarrow}$ & $\longrightarrow$ & $\underset{Aut}{\nrightarrow}$ & $\underset{Unital}{\nrightarrow}$ & $\underset{Aut}{\nrightarrow}$ & $\underset{Aut}{\nrightarrow}$ & $\underset{trL_x}{\nrightarrow}$ & $\longrightarrow$ & $\underset{Aut}{\nrightarrow}$ & $\underset{Aut}{\nrightarrow}$ & $\underset{Aut}{\nrightarrow}$ & $\underset{Ass}{\nrightarrow}$ & $\underset{Ass}{\nrightarrow}$ & $\underset{Aut}{\nrightarrow}$ & $\underset{Ass}{\nrightarrow}$ & $\underset{Aut}{\nrightarrow}$ & $\underset{t}{\longrightarrow}$ & $\underset{Ass}{\nrightarrow}$ & $\underset{Ass}{\nrightarrow}$ & $\underset{Ass}{\nrightarrow}$ & $\underset{Aut}{\nrightarrow}$ & $\underset{Aut}{\nrightarrow}$ & $\underset{Aut}{\nrightarrow}$ & $\underset{Aut}{\nrightarrow}$  & $\underset{Ass}{\nrightarrow}$ &  $\underset{Ass}{\nrightarrow}$\tabularnewline
					\hline 
					$J_{3,11}^{0}$ & $\underset{trivial}{\longrightarrow}$ & $\underset{t}{\longrightarrow}$ & $\underset{t}{\longrightarrow}$ & $\underset{t}{\longrightarrow}$ & $\underset{Aut}{\nrightarrow}$ & $\underset{Unital}{\nrightarrow}$ & $\underset{Unital}{\nrightarrow}$ & $\longrightarrow$ & $\underset{t}{\longrightarrow}$ & $\longrightarrow$ & $\longrightarrow$ & $\underset{Aut}{\nrightarrow}$ & $\underset{Ass}{\nrightarrow}$ & $\underset{Ass}{\nrightarrow}$ & $\underset{Ass}{\nrightarrow}$ & $\underset{Aut}{\nrightarrow}$ & $\underset{Ass}{\nrightarrow}$ & $\underset{Ass}{\nrightarrow}$ & $\underset{t}{\longrightarrow}$ & $\underset{Ass}{\nrightarrow}$ & $\underset{Ass}{\nrightarrow}$ & $\underset{Ass}{\nrightarrow}$ & $\underset{Aut}{\nrightarrow}$ & $\underset{Ass}{\nrightarrow}$ & $\underset{Aut}{\nrightarrow}$ & $\underset{Aut}{\nrightarrow}$  & $\underset{Ass}{\nrightarrow}$ &  $\underset{Ass}{\nrightarrow}$\tabularnewline
					\hline 
					$J_{3,13}^{0}$ & $\underset{trivial}{\longrightarrow}$ & $\underset{t}{\longrightarrow}$ & $\underset{t}{\longrightarrow}$ & $\underset{t}{\longrightarrow}$ & $\longrightarrow$ & $\underset{t}{\longrightarrow}$ & $\underset{t}{\longrightarrow}$ & $\underset{t}{\longrightarrow}$ & $\underset{t}{\longrightarrow}$ & $\underset{t}{\longrightarrow}$ & $\longrightarrow$ & $\longrightarrow$ & $\underset{Ass}{\nrightarrow}$ & $\underset{Ass}{\nrightarrow}$ & $\underset{Ass}{\nrightarrow}$ & $\underset{Ass}{\nrightarrow}$ & $\underset{Ass}{\nrightarrow}$ & $\underset{Ass}{\nrightarrow}$ & $\underset{t}{\longrightarrow}$ & $\underset{Ass}{\nrightarrow}$ & $\underset{Ass}{\nrightarrow}$ & $\underset{Ass}{\nrightarrow}$ & $\underset{Ass}{\nrightarrow}$ & $\underset{Ass}{\nrightarrow}$ & $\underset{Ass}{\nrightarrow}$ & $\underset{Ass}{\nrightarrow}$  & $\underset{Ass}{\nrightarrow}$ &  $\underset{Ass}{\nrightarrow}$\tabularnewline
					\hline 
					$J_{3,15}^1$ & $\underset{trivial}{\longrightarrow}$ & $\longrightarrow$ & $\longrightarrow$ & $\underset{x^2y=0}{\nrightarrow}$ & $\underset{Aut}{\nrightarrow}$ & $\underset{Unital}{\nrightarrow}$ & $\underset{Aut}{\nrightarrow}$ & $\underset{Aut}{\nrightarrow}$ & $\underset{x^2y=0}{\nrightarrow}$ & $\underset{Aut}{\nrightarrow}$ & $\underset{Aut}{\nrightarrow}$ & $\underset{Aut}{\nrightarrow}$ & $\longrightarrow$ & $\longrightarrow$ & $\underset{x^2y=0}{\nrightarrow}$ & $\underset{Aut}{\nrightarrow}$ & $\underset{x^2y=0}{\nrightarrow}$ & $\underset{Aut}{\nrightarrow}$ & $\longrightarrow$ & $\longrightarrow$ & $\underset{dimZ}{\nrightarrow}$ & $\underset{dimZ}{\nrightarrow}$ & $\underset{Aut}{\nrightarrow}$ & $\underset{Aut}{\nrightarrow}$ & $\underset{Aut}{\nrightarrow}$ & $\underset{Aut}{\nrightarrow}$  & $\underset{Ann}{\nrightarrow}$ &  $\underset{Ann}{\nrightarrow}$\tabularnewline
					\hline 
					$J_{3,15}^0$ & $\underset{trivial}{\longrightarrow}$ & $\longrightarrow$ & $\longrightarrow$ & $\underset{Aut}{\nrightarrow}$ & $\underset{Aut}{\nrightarrow}$ & $\underset{Unital}{\nrightarrow}$ & $\underset{Aut}{\nrightarrow}$ & $\underset{Aut}{\nrightarrow}$ & $\underset{x^2y=0}{\nrightarrow}$ & $\underset{Aut}{\nrightarrow}$ & $\underset{Aut}{\nrightarrow}$ & $\underset{Aut}{\nrightarrow}$ & $\underset{Aut}{\nrightarrow}$ & $\longrightarrow$ & $\underset{Unital}{\nrightarrow}$ & $\underset{Aut}{\nrightarrow}$ & $\underset{Aut}{\nrightarrow}$ & $\underset{Aut}{\nrightarrow}$ & $\longrightarrow$ & $\longrightarrow$ & $\underset{Ann}{\nrightarrow}$ & $\underset{Ann}{\nrightarrow}$ & $\underset{Aut}{\nrightarrow}$ & $\underset{Aut}{\nrightarrow}$ & $\underset{Aut}{\nrightarrow}$ & $\underset{Aut}{\nrightarrow}$  & $\underset{Ann}{\nrightarrow}$ &  $\underset{Ann}{\nrightarrow}$\tabularnewline
					\hline 
					$J_{3,16}^{1,0}$ & $\underset{trivial}{\longrightarrow}$ & $\longrightarrow$ & $\longrightarrow$ & $\underset{Alicyclic}{\nrightarrow}$ & $\underset{Aut}{\nrightarrow}$ & $\longrightarrow$ & $\underset{Aut}{\nrightarrow}$ & $\underset{Aut}{\nrightarrow}$ & $\underset{L_{x^2}=c(x)I}{\nrightarrow}$ & $\underset{Aut}{\nrightarrow}$ & $\underset{Aut}{\nrightarrow}$ & $\underset{Aut}{\nrightarrow}$ & $\underset{Aut}{\nrightarrow}$ &  $\longrightarrow$  & $\longrightarrow$ & $\underset{Aut}{\nrightarrow}$ & $\underset{tr(L_{xy})}{\nrightarrow}$ & $\underset{Aut}{\nrightarrow}$ & $\longrightarrow$ & $\longrightarrow$ & $\underset{dimZ}{\nrightarrow}$ & $\underset{dimZ}{\nrightarrow}$ & $\underset{Aut}{\nrightarrow}$ & $\underset{Aut}{\nrightarrow}$ & $\underset{Aut}{\nrightarrow}$ & $\underset{Aut}{\nrightarrow}$  & $\underset{dim Z}{\nrightarrow}$ &  $\underset{dim Z}{\nrightarrow}$\tabularnewline
					\hline 
					$J_{3,16}^{1,1}$ & $\underset{trivial}{\longrightarrow}$ & $\longrightarrow$ & $\underset{t}{\longrightarrow}$ & $\underset{Alicyclic}{\nrightarrow}$ & $\underset{Aut}{\nrightarrow}$ & $\longrightarrow$ & $\underset{Alicyclic}{\nrightarrow}$ & $\underset{L_{x^2}=c(x)I}{\nrightarrow}$ & $\underset{L_{x^2}=c(x)I}{\nrightarrow}$ & $\underset{L_{x^2}=c(x)I}{\nrightarrow}$ & $\underset{Aut}{\nrightarrow}$ & $\underset{Aut}{\nrightarrow}$ & $\longrightarrow$ & $\underset{t}{\longrightarrow}$ & $\longrightarrow$ & $\longrightarrow$ & $\longrightarrow$ & $\underset{L_{x^2}=c(x)I}{\nrightarrow}$ & $\longrightarrow$ & $\underset{t}{\longrightarrow}$ & $\underset{dimZ}{\nrightarrow}$ & $\underset{dimZ}{\nrightarrow}$ & $\underset{Aut}{\nrightarrow}$ & $\underset{dimZ}{\nrightarrow}$ & $\underset{Aut}{\nrightarrow}$ & $\underset{Aut}{\nrightarrow}$ & $\underset{dim Z}{\nrightarrow}$ &  $\underset{dim Z}{\nrightarrow}$\tabularnewline
					\hline 
					$J_{3,17}$ & $\underset{trivial}{\longrightarrow}$ & $\underset{t}{\longrightarrow}$ & $\longrightarrow$ & $\underset{Aut}{\nrightarrow}$ & $\underset{Aut}{\nrightarrow}$ & $\longrightarrow$ & $\underset{Aut}{\nrightarrow}$ & $\underset{Aut}{\nrightarrow}$ & $\underset{L_{x^2}=c(x)I}{\nrightarrow}$ & $\underset{Aut}{\nrightarrow}$ & $\underset{Aut}{\nrightarrow}$ & $\underset{Aut}{\nrightarrow}$ & $\underset{Aut}{\nrightarrow}$ & $\underset{Aut}{\nrightarrow}$ & $\underset{Characteristic}{\nrightarrow}$ & $\underset{Aut}{\nrightarrow}$ & $\longrightarrow$ & $\underset{Aut}{\nrightarrow}$ & $\longrightarrow$ & $\underset{Characteristic}{\nrightarrow}$ & $\underset{dimZ}{\nrightarrow}$ & $\underset{dimZ}{\nrightarrow}$ & $\underset{Aut}{\nrightarrow}$ & $\underset{Aut}{\nrightarrow}$ & $\underset{Aut}{\nrightarrow}$ & $\underset{Aut}{\nrightarrow}$ & $\underset{dim Z}{\nrightarrow}$ &  $\underset{dim Z}{\nrightarrow}$\tabularnewline
					\hline 
					$J_{3,24}$ & $\underset{trivial}{\longrightarrow}$ & $\underset{t}{\longrightarrow}$ & $\underset{t}{\longrightarrow}$ & $\underset{Alicyclic}{\nrightarrow}$ & $\underset{Aut}{\nrightarrow}$ & $\underset{Unital}{\nrightarrow}$ & $\underset{Aut}{\nrightarrow}$ & $\underset{Aut}{\nrightarrow}$ & $\longrightarrow$ & $\underset{Aut}{\nrightarrow}$ & $\underset{Aut}{\nrightarrow}$ & $\underset{Aut}{\nrightarrow}$ & $\underset{Aut}{\nrightarrow}$ & $\longrightarrow$ & $\underset{Unital}{\nrightarrow}$ & $\underset{Unital}{\nrightarrow}$ & $\underset{Unital}{\nrightarrow}$ & $\longrightarrow$ & $\longrightarrow$ & $\longrightarrow$ & $\underset{dimZ}{\nrightarrow}$ & $\underset{dimZ}{\nrightarrow}$ & $\underset{Aut}{\nrightarrow}$ & $\underset{Aut}{\nrightarrow}$ & $\underset{Aut}{\nrightarrow}$ & $\underset{Aut}{\nrightarrow}$ & $\underset{dim Z}{\nrightarrow}$ &  $\underset{dim A^2}{\nrightarrow}$\tabularnewline
					\hline 
					$J_{3,25}$ & $\underset{trivial}{\longrightarrow}$ & $\underset{Aut}{\nrightarrow}$ & $\underset{Aut}{\nrightarrow}$ & $\underset{Aut}{\nrightarrow}$ & $\underset{Aut}{\nrightarrow}$ & $\underset{Aut}{\nrightarrow}$ & $\underset{Aut}{\nrightarrow}$ & $\underset{Aut}{\nrightarrow}$ & $\underset{Aut}{\nrightarrow}$ & $\underset{Aut}{\nrightarrow}$ & $\underset{Aut}{\nrightarrow}$ & $\underset{Aut}{\nrightarrow}$ & $\underset{Aut}{\nrightarrow}$ & $\underset{Aut}{\nrightarrow}$ & $\underset{x^2=0}{\nrightarrow}$ & $\underset{Aut}{\nrightarrow}$ & $\underset{Aut}{\nrightarrow}$ & $\underset{Aut}{\nrightarrow}$ & $\longrightarrow$ & $\underset{Aut}{\nrightarrow}$ & $\underset{Aut}{\nrightarrow}$ & $\underset{Aut}{\nrightarrow}$ & $\underset{Aut}{\nrightarrow}$ & $\underset{Aut}{\nrightarrow}$ & $\underset{Aut}{\nrightarrow}$ & $\underset{Aut}{\nrightarrow}$ & $\underset{Ann}{\nrightarrow}$ &  $\underset{Ann}{\nrightarrow}$\tabularnewline
					\hline 
					$J_{3,26}$ & $\underset{trivial}{\longrightarrow}$ & $\underset{x{{}^2}=0}{\nrightarrow}$ & $\underset{Aut}{\nrightarrow}$ & $\underset{Aut}{\nrightarrow}$ & $\underset{Aut}{\nrightarrow}$ & $\underset{Aut}{\nrightarrow}$ & $\underset{Aut}{\nrightarrow}$ & $\underset{Aut}{\nrightarrow}$ & $\underset{Aut}{\nrightarrow}$ & $\underset{Aut}{\nrightarrow}$ & $\underset{Aut}{\nrightarrow}$ & $\underset{Aut}{\nrightarrow}$ & $\underset{Aut}{\nrightarrow}$ & $\underset{Aut}{\nrightarrow}$ & $\underset{Aut}{\nrightarrow}$ & $\underset{Aut}{\nrightarrow}$ & $\underset{Aut}{\nrightarrow}$ & $\underset{Aut}{\nrightarrow}$ & $\longrightarrow$ & $\longrightarrow$ & $\underset{dimZ}{\nrightarrow}$ & $\underset{Aut}{\nrightarrow}$ & $\underset{Aut}{\nrightarrow}$ & $\underset{Aut}{\nrightarrow}$ & $\underset{Aut}{\nrightarrow}$ & $\underset{Aut}{\nrightarrow}$ & $\underset{Ann}{\nrightarrow}$ &  $\underset{Ann}{\nrightarrow}$\tabularnewline
					\hline 
					$J_{3,27}$ & $\underset{trivial}{\longrightarrow}$ & $\underset{Aut}{\nrightarrow}$ & $\underset{Aut}{\nrightarrow}$ & $\underset{Aut}{\nrightarrow}$ & $\underset{Aut}{\nrightarrow}$ & $\underset{Aut}{\nrightarrow}$ & $\underset{Aut}{\nrightarrow}$ & $\underset{Aut}{\nrightarrow}$ & $\underset{Aut}{\nrightarrow}$ & $\underset{Aut}{\nrightarrow}$ & $\underset{Aut}{\nrightarrow}$ & $\underset{Aut}{\nrightarrow}$ & $\underset{Aut}{\nrightarrow}$ & $\underset{Aut}{\nrightarrow}$ & $\underset{Aut}{\nrightarrow}$ & $\underset{Aut}{\nrightarrow}$ & $\underset{Aut}{\nrightarrow}$ & $\underset{Aut}{\nrightarrow}$ & $\underset{Aut}{\nrightarrow}$ & $\underset{Aut}{\nrightarrow}$ & $\longrightarrow$ & $\underset{Aut}{\nrightarrow}$ & $\underset{Aut}{\nrightarrow}$ & $\underset{Aut}{\nrightarrow}$ & $\underset{Aut}{\nrightarrow}$ & $\underset{Aut}{\nrightarrow}$ & $\underset{Aut}{\nrightarrow}$ &  $\underset{dim A^2}{\nrightarrow}$\tabularnewline
					\hline 
					$J_{3,28}^{\gamma\neq0}$ & $\underset{trivial}{\longrightarrow}$ & $\underset{x{{}^2}=0}{\nrightarrow}$ & $\underset{Aut}{\nrightarrow}$ & $\underset{Aut}{\nrightarrow}$ & $\underset{Aut}{\nrightarrow}$ & $\underset{Aut}{\nrightarrow}$ & $\underset{Aut}{\nrightarrow}$ & $\underset{Aut}{\nrightarrow}$ & $\underset{Aut}{\nrightarrow}$ & $\underset{Aut}{\nrightarrow}$ & $\underset{Aut}{\nrightarrow}$ & $\underset{Aut}{\nrightarrow}$ & $\underset{Aut}{\nrightarrow}$ & $\underset{Aut}{\nrightarrow}$ & $\underset{Aut}{\nrightarrow}$ & $\underset{Aut}{\nrightarrow}$ & $\underset{Aut}{\nrightarrow}$ & $\underset{Aut}{\nrightarrow}$ & $\longrightarrow$ & $\underset{Aut}{\nrightarrow}$ & $\underset{Characteristic}{\nrightarrow}$ & $\longrightarrow$ & $\underset{Aut}{\nrightarrow}$ & $\underset{Aut}{\nrightarrow}$ & $\underset{Aut}{\nrightarrow}$ & $\underset{Aut}{\nrightarrow}$ & $\underset{Characteristic}{\nrightarrow}$ &  $\underset{dim A^2}{\nrightarrow}$\tabularnewline
					\hline
					$J_{3,29}$ & $\underset{trivial}{\longrightarrow}$ & $\underset{x{{}^2}=0}{\nrightarrow}$ & $\underset{x{{}^2}=0}{\nrightarrow}$ & $\underset{x{{}^2}=0}{\nrightarrow}$ & $\underset{Aut}{\nrightarrow}$ & $\underset{x{{}^2}=0}{\nrightarrow}$ & $\underset{x{{}^2}=0}{\nrightarrow}$ & $\underset{x{{}^2}=0}{\nrightarrow}$ & $\underset{x{{}^2}=0}{\nrightarrow}$ & $\underset{x{{}^2}=0}{\nrightarrow}$ & $\underset{Aut}{\nrightarrow}$ & $\underset{Aut}{\nrightarrow}$ & $\underset{x{{}^2}=0}{\nrightarrow}$ & $\underset{x{{}^2}=0}{\nrightarrow}$ & $\underset{x{{}^2}=0}{\nrightarrow}$ & $\underset{Aut}{\nrightarrow}$ & $\underset{x{{}^2}=0}{\nrightarrow}$ & $\underset{x{{}^2}=0}{\nrightarrow}$ & $\longrightarrow$ & $\longrightarrow$ & $\longrightarrow$ & $\longrightarrow$ & $\longrightarrow$ & $\longrightarrow$ & $\underset{Aut}{\nrightarrow}$ & $\underset{Aut}{\nrightarrow}$ & ${\longrightarrow}$ &  $\underset{dim A^2}{\nrightarrow}$\tabularnewline
					\hline 
					$J_{3,30}$ & $\underset{trivial}{\longrightarrow}$ & $\underset{x{{}^2}=0}{\nrightarrow}$ & $\underset{x{{}^2}=0}{\nrightarrow}$ & $\underset{x{{}^2}=0}{\nrightarrow}$ & $\underset{Aut}{\nrightarrow}$ & $\underset{x{{}^2}=0}{\nrightarrow}$ & $\underset{Aut}{\nrightarrow}$ & $\underset{Aut}{\nrightarrow}$ & $\underset{x{{}^2}=0}{\nrightarrow}$ & $\underset{Aut}{\nrightarrow}$ & $\underset{Aut}{\nrightarrow}$ & $\underset{Aut}{\nrightarrow}$ & $\underset{Aut}{\nrightarrow}$ & $\underset{x{{}^2}=0}{\nrightarrow}$ & $\underset{x{{}^2}=0}{\nrightarrow}$ & $\underset{Aut}{\nrightarrow}$ & $\underset{x{{}^2}=0}{\nrightarrow}$ & $\underset{Aut}{\nrightarrow}$ & $\longrightarrow$ & $\longrightarrow$ & $\longrightarrow$ & $\longrightarrow$ & $\underset{Aut}{\nrightarrow}$ & $\longrightarrow$ & $\underset{Aut}{\nrightarrow}$ & $\underset{Aut}{\nrightarrow}$ & ${\longrightarrow}$ &  $\underset{dim A^2}{\nrightarrow}$\tabularnewline
					\hline 
					$J_{3,31}^{\gamma\neq0}$ & $\underset{trivial}{\longrightarrow}$ & $\underset{x{{}^2}=0}{\nrightarrow}$ & $\underset{x{{}^2}=0}{\nrightarrow}$ & $\underset{x{{}^2}=0}{\nrightarrow}$ & $\underset{Aut}{\nrightarrow}$ & $\underset{x{{}^2}=0}{\nrightarrow}$ & $\underset{x{{}^2}=0}{\nrightarrow}$ & $\underset{x{{}^2}=0}{\nrightarrow}$ & $\underset{x{{}^2}=0}{\nrightarrow}$ & $\underset{x{{}^2}=0}{\nrightarrow}$ & $\underset{Aut}{\nrightarrow}$ & $\underset{Aut}{\nrightarrow}$ & $\underset{x{{}^2}=0}{\nrightarrow}$ & $\underset{x{{}^2}=0}{\nrightarrow}$ & $\underset{x{{}^2}=0}{\nrightarrow}$ & $\underset{Aut}{\nrightarrow}$ & $\underset{x{{}^2}=0}{\nrightarrow}$ & $\underset{x{{}^2}=0}{\nrightarrow}$ & $\longrightarrow$ & $\longrightarrow$ & $\longrightarrow$ & $\longrightarrow$ & $\underset{Aut}{\nrightarrow}$ & $\longrightarrow$ & $\longrightarrow$ & $\underset{Aut}{\nrightarrow}$ & ${\longrightarrow}$ &  $\underset{Det}{\nrightarrow}$\tabularnewline
					\hline 
					$J_{3,32}$ & $\underset{trivial}{\longrightarrow}$ & $\underset{x{{}^2}=0}{\nrightarrow}$ & $\underset{x{{}^2}=0}{\nrightarrow}$ & $\underset{x{{}^2}=0}{\nrightarrow}$ & $\underset{Aut}{\nrightarrow}$ & $\underset{x{{}^2}=0}{\nrightarrow}$ & $\underset{x{{}^2}=0}{\nrightarrow}$ & $\underset{x{{}^2}=0}{\nrightarrow}$ & $\underset{x{{}^2}=0}{\nrightarrow}$ & $\underset{x{{}^2}=0}{\nrightarrow}$ & $\underset{Aut}{\nrightarrow}$ & $\underset{Aut}{\nrightarrow}$ & $\underset{x{{}^2}=0}{\nrightarrow}$ & $\underset{x{{}^2}=0}{\nrightarrow}$ & $\underset{x{{}^2}=0}{\nrightarrow}$ & $\underset{Aut}{\nrightarrow}$ & $\underset{x{{}^2}=0}{\nrightarrow}$ & $\underset{x{{}^2}=0}{\nrightarrow}$ & $\longrightarrow$ & $\longrightarrow$ & $\longrightarrow$ & $\longrightarrow$ & $\underset{Aut}{\nrightarrow}$ & $\longrightarrow$ & $\underset{Aut}{\nrightarrow}$ & $\longrightarrow$ & ${\longrightarrow}$ &  ${\longrightarrow}$\tabularnewline
					\hline 
					$J_{3,28}^{0}$ & $\underset{trivial}{\longrightarrow}$ & $\underset{x{{}^2}=0}{\nrightarrow}$  & $\underset{x{{}^2}=0}{\nrightarrow}$  & $\underset{x{{}^2}=0}{\nrightarrow}$  & $\underset{x{{}^2}=0}{\nrightarrow}$  & $\underset{x{{}^2}=0}{\nrightarrow}$  & $\underset{x{{}^2}=0}{\nrightarrow}$  & $\underset{x{{}^2}=0}{\nrightarrow}$  & $\underset{x{{}^2}=0}{\nrightarrow}$  & $\underset{x{{}^2}=0}{\nrightarrow}$  & $\underset{x{{}^2}=0}{\nrightarrow}$  & $\underset{x{{}^2}=0}{\nrightarrow}$  & $\underset{x{{}^2}=0}{\nrightarrow}$  & $\underset{x{{}^2}=0}{\nrightarrow}$  & $\underset{x{{}^2}=0}{\nrightarrow}$  & $\underset{x{{}^2}=0}{\nrightarrow}$  & $\underset{x{{}^2}=0}{\nrightarrow}$  & $\underset{x{{}^2}=0}{\nrightarrow}$  & ${\longrightarrow}$ & $\underset{trL_x}{\nrightarrow}$ & ${\longrightarrow}$ & $\underset{trL_x}{\nrightarrow}$ & $\underset{trL_x}{\nrightarrow}$ & $\underset{trL_x}{\nrightarrow}$ & $\underset{trL_x}{\nrightarrow}$ & $\underset{trL_x}{\nrightarrow}$ & ${\longrightarrow}$ &  $\underset{dim A^2}{\nrightarrow}$\tabularnewline
					\hline 
					$J_{3,31}^{0}$ & $\underset{trivial}{\longrightarrow}$  & $\underset{x{{}^2}=0}{\nrightarrow}$  & $\underset{x{{}^2}=0}{\nrightarrow}$  & $\underset{x{{}^2}=0}{\nrightarrow}$  & $\underset{x{{}^2}=0}{\nrightarrow}$  & $\underset{x{{}^2}=0}{\nrightarrow}$  & $\underset{x{{}^2}=0}{\nrightarrow}$  & $\underset{x{{}^2}=0}{\nrightarrow}$  & $\underset{x{{}^2}=0}{\nrightarrow}$  & $\underset{x{{}^2}=0}{\nrightarrow}$  & $\underset{x{{}^2}=0}{\nrightarrow}$  & $\underset{x{{}^2}=0}{\nrightarrow}$  & $\underset{x{{}^2}=0}{\nrightarrow}$  & $\underset{x{{}^2}=0}{\nrightarrow}$  & $\underset{x{{}^2}=0}{\nrightarrow}$  & $\underset{x{{}^2}=0}{\nrightarrow}$  & $\underset{x{{}^2}=0}{\nrightarrow}$  & $\underset{x{{}^2}=0}{\nrightarrow}$  & ${\longrightarrow}$ & $\underset{trL_x}{\nrightarrow}$ & ${\longrightarrow}$ & $\underset{trL_x}{\nrightarrow}$ & $\underset{trL_x}{\nrightarrow}$ & $\underset{trL_x}{\nrightarrow}$ & $\underset{trL_x}{\nrightarrow}$ & $\underset{trL_x}{\nrightarrow}$ & ${\longrightarrow}$ &  ${\longrightarrow}$\tabularnewline
					\hline 
					
				\end{tabular}
		\end{adjustbox} }
		\caption{Deformations of three-dimensional Jordan algebras.}
		\label{classification}
	\end{table}
	
	\clearpage
	The deformations used in the Table \ref{classification} are given in Tables \ref{table1} and \ref{table2} below.
	
	\begin{table}[!h]
		\resizebox{1\columnwidth}{!}{%
			\begin{tabular}{|l|l|}
				\hline
				Deformation & Change of basis \\ \hline
				$J_{3,3}^{1}\to J_{3,2}$   &     $f_1=ta,\ f_2=t^2c,\ f_3=t^2b$              \\ \hline
				$J_{3,6}\to J_{3,2}$   &       $f_1=t+a,\ f_2=t^2,\ f_3=b$            \\ \hline
				$J_{3,9} \to J_{3,2}$   &        $f_1=ta+b,\ f_2=t^2a,\ f_3=c$           \\ \hline
				$J_{3,4} \to J_{3,2}$   &        $f_1=ta,\ f_2=t^2b,\ f_3=c$           \\ \hline
				$J_{3,4} \to J_{3,3}^1$   &       $f_1=ta+tc,\ f_2=ta+tb ,\ f_3=t^2b  $       \\ \hline
				$J_{3,15}^0 \to J_{3,25}$   &        $f_1=ta+a^2,\ f_2=tb+a^2,\ f_3=t^2b$           \\ \hline
				$J_{3,15}^0 \to J_{3,2}$   &        $f_1=ta,\ f_2=t^2a^2,\ f_3=tb$           \\ \hline
				$J_{3,15}^0 \to J_{3,3}^1$   &       $f_1=ta+ta^2, f_2=ta+tb,f_3=t^2a$          \\ \hline
				$J_{3,16}^{1,1} \to J_{3,25}$   &        $f_1=t+tb,\ f_2=a+b,\ f_3=t$           \\ \hline
				$J_{3,16}^{1,1} \to J_{3,2}$   &        $f_1=t+a+b,\ f_2=t^2,\ f_3=t^3a$           \\ \hline
				$J_{3,16}^{1,1} \to J_{3,6}$   &        $f_1=1,\ f_2=ta,\ f_3=tb$           \\ \hline
				$J_{3,17} \to J_{3,25}$   &        $f_1=t+a,\ f_2=tb,\ f_3=t$           \\ \hline
				$J_{3,17} \to J_{3,3}^1$   &        $f_1=t+ta,\ f_2=t+tb,\ f_3=t^2$           \\ \hline
				$J_{3,17} \to J_{3,6}$   &        $f_1=1,\ f_2=ta,\ f_3=tb$           \\ \hline
				$J_{3,7} \to J_{3,4}$   &   
				$f_1=a+t,\ f_2=a^2+t^2,\ f_3=t^3+t^2a+ta^2$           \\ \hline
				$J_{3,7} \to J_{3,6}$   &        $f_1=1,\ f_2=ta,\ f_3=a^2$           \\ \hline
				$J_{3,8} \to J_{3,4}$   &      $f_1=a+tc,\ f_2=b+t^2c,\ f_3=t^3c$           \\ \hline
				$J_{3,8} \to J_{3,9}$   &      $f_1=c,\ f_2=ta,\ f_3=b$ \            \\ \hline
				$J_{3,10}^0 \to J_{3,4}$   &      $f_1=ta+b+c,\ f_2=t^2a,\ f_3=t^3a+t^2b$            \\ \hline
				$J_{3,15}^1 \to J_{3,25}$   &      $f_1=tb+(t+1)a^2,\ f_2=t(a+b+a^2),\ f_3=t^2(b+a^2)$   \\ \hline
				$J_{3,15}^1 \to J_{3,15}^0$   &    $f_1=a+(t+1)b+(t+1)a^2,\ f_2=t^2a^2,\ f_3=t^2(b+a^2)$    \\ \hline
				$J_{3,15}^1 \to J_{3,26}$   &      $f_1=tb+(t+1)a^2,\ f_2=a+b,\ f_3=t(b+a^2)$     \\ \hline
				$J_{3,15}^1 \to J_{3,2}$   &      $f_1=b+a^2,\ f_2=a^2,\ f_3=t(a+b+a^2)$   \\ \hline
				$J_{3,24} \to J_{3,9}$   &      $f_1=a,\ f_2=b,\ f_3=tc$  \          \\ \hline
				$J_{3,5} \to J_{3,8}$   &      $f_1=t+a,\ f_2=t^2,\ f_3=b$    \        \\ \hline
				$J_{3,5} \to J_{3,7}$   &      $f_1=1,\ f_2=a+tb,\ f_3=t^2b$    \        \\ \hline
				$J^0_{3,11} \to J_{3,8}$   &      $f_1=ta+c,\ f_2=t^2a,\ f_3=b$ \           \\ \hline
				$J^0_{3,11} \to J^0_{3,10}$   &      $f_1=a,\ f_2=t^2b,\ f_3=c$  \          \\ \hline
				$J^0_{3,13} \to J_{3,5}$   &      $f_1=a+c,\  f_2=tb+tc,\ f_3=b$           \\ \hline
				$J^0_{3,13} \to J^0_{3,11}$   &      $f_1=a,\ f_2=b,\ f_3=tc$           \\ \hline
				$J^{1,1}_{3,16} \to J^{1,0}_{3,16}$   &      $f_1=1,\ f_2=a,\ f_3=tb$          \\ \hline
				$J^{1,1}_{3,16} \to J_{3,15}^1$   & 
				$f_1=1+ta+b,\ f_2=t^2,\ f_3=ta+t^2$          \\ \hline
				$J^{1,1}_{3,16} \to J_{3,17}$   & 
				$f_1=1,f_2=ta,\ f_3=t^{-1}(a+b)$         \\ \hline
				$J^{1,0}_{3,16} \to J_{3,2}$   & 
				$f_1=t+b,f_2=t^2,\ f_3=t^3a$           \\ \hline
				$J_{3,32} \to J_{3,25}$   & 
				$f_1=t^2a,f_2=tb,\ f_3=t^3c$           \\ \hline
				$J^{\gamma\neq0}_{3,31} \to J_{3,26}$   & 
				$f_1=t(\gamma a+tb),f_2=a+\gamma^{-1}(b+c),\ f_3=t(b+c)$           \\ \hline
				$J_{3,32} \to J_{3,30}$   & 
				$f_1=xa+c,f_2=(1+x+t)a+(1+x)b+xc,\ f_3=(1+x)a+(1+x)b+xc$ for $x\neq 1, x^3=1$.           \\ \hline
				$J^{\gamma\neq0}_{3,31} \to J_{3,30}$   & 
				$f_1=a+\gamma^{-1}b+(\gamma^{-1}+t)c,f_2=\gamma^{-1}(b+c),\ f_3=\gamma^{-1}b+(\gamma^{-1}+t)c$       \\ \hline
				$J_{3,5} \to J_{3,10}^0$   & 
				$f_1=1+b,f_2=a,f_3=t$  \      \\ \hline
				$J_{3,28}^0 \to J_{3,25}$ &
				$f_1=ta,\ f_2=tb,\ f_3=t^2c$ \\ \hline
				$J_{3,31}^0 \to J_{3,25}$ &
				$f_1=ta,\ f_2=tb,\ f_3=t^2c$ \\ \hline
				$J_{3,28}^0 \to J_{3,27}$ &
				$f_1=b+c,\ f_2=tb,\ f_3=a$ \\ \hline
				$J_{3,31}^0 \to J_{3,27}$ &
				$f_1=b+c,\ f_2=tb,\ f_3=a$ \\ \hline
				$J_{3,31}^0 \to J_{3,28}^0$ &
				$f_1=a,\ f_2=tb,\ f_3=tc$ \\ \hline
				$J_{3,29} \to J_{3,28}^0$ &
				$f_1=a+b,\ f_2=tc,\ f_3=ta$ \\ \hline
				$J_{3,30} \to J_{3,28}^0$ &
				$f_1=a+b,\ f_2=tc,\ f_3=ta$ \\ \hline
				$J_{3,31}^{\gamma\neq 0} \to J_{3,28}^0$ &
				$f_1=a,\ f_2=tb,\ f_3=tc$ \\ \hline
				$J_{3,31}^0 \to J_{3,28}^0$ &
				$f_1=a,\ f_2=tb,\ f_3=tc$ \\ \hline
				$J_{3,32} \to J_{3,28}^0$ &
				$f_1=a,\ f_2=tb,\ f_3=tc$ \\ \hline
				$J_{3,32} \to J_{3,31}^0$ &
				$f_1=a+t^{-1}(b+c),\
				f_2=a+tc,\
				f_3=a+tc+t^{-1}(b+c)$ \\ \hline       
			\end{tabular}
		}
		\caption{Existence of deformations in $Jor_3$}\label{table1}
	\end{table}
	
	\clearpage

	\begin{table}[h!]
		\begin{adjustbox}{max totalsize={\textwidth}{1\textheight},center}
			\begin{tabular}{|l|l|}
				\hline
				Deformation & Change of basis \\ \hline
				$J_{3,7} \to J_{3,2}$   & 
				$f_1=a,f_2=a^2,f_3=t$ \            \\ \hline
				$J^0_{3,11} \to J_{3,8}$   & 
				$f_1=ta+c,f_2=t^2a,f_3=a+b$  \           \\ \hline
				$J_{3,7} \to J^1_{3,3}$   & 
				$f_1=t^{-1}a,f_2=t+t^{-3}a^2,f_3=t^{-2}a^2$     \\ \hline
				$J_{3,4} \to J_{3,25}$   & 
				$f_1=ta,f_2=b,f_3=tc$\            \\ \hline
				$J^1_{3,3} \to J_{3,25}$   & 
				$f_1=ta,f_2=a+b,f_3=tc$ \         \\ \hline
				$J^{1,0}_{3,16} \to J_{3,25}$   & 
				$f_1=t+ta+b,f_2=t^ 2a,f_3=t^3$ \              \\ \hline
				$J_{3,15}^0 \to J_{3,26}$           &  $f_1=t^{-1}a^2+b,f_2=a,f_3=b$      
				\\ \hline
				$J_{3,24} \to J^0_{3,15}$          & $f_1=ta+c,\,f_2=t^2a,\,f_3=b$ \       \\  \hline  
				$J_{3,24} \to J_{3,25}$          & $f_1=tc,\,f_2=t^2a+b,\,f_3=tb$  \                                     \\ \hline
				$J_{3,24} \to J_{3,26}$          &    $f_1=ta+b,f_2=c,f_3=b$                                   \\ \hline
				$J_{3,26} \to J_{3,25}$           &      $f_1=a,f_2=tb,f_3=tc$    \                               \\ \hline
				$J_{3,28}^{\gamma\neq0} \to J_{3,25}$ &  $f_1=b,f_2=t a,f_3=t c$   \                   \\ \hline
				$J_{3,29} \to  J_{3,25}$ & $f_1=a,f_2=tb,f_3=tc$           
				\\ \hline
				$J_{3,29} \to   J_{3,26}$        &  $f_1=b+ta,f_2=c,f_3=b$\                                  \\ \hline
				$J_{3,29} \to   J_{3,27}$        & $f_1=tb+tc,f_2=t^3b,f_3=a+b$   \   \\ \hline
				$J_{3,29} \to  J_{3,28}^{\gamma\neq0}$       &   $f_1=a+b+\gamma c,f_2=t c,f_3=t b$     \                            \\ \hline
				$J_{3,29} \to  J_{3,30}$         & $f_1=a+c,f_2=(t+1)b+c,f_3=b+c$ \   \\ \hline
				$J_{3,30} \to  J_{3,25}$        & $f_1=ta,f_2=tb,f_3=t^2c$ \                                   \\ \hline
				$J_{3,30} \to  J_{3,26}$         & $f_1=t^2b+tc,f_2=a,f_3=tc$   \                            \\ \hline
				$J_{3,30} \to  J_{3,27}$         &  $f_1=b+c,f_2=tc,f_3=a+c$   \                              \\ \hline
				$J_{3,30} \to  J_{3,28}^{\gamma\neq0}$       &    $f_1=\lambda a+\mu c,f_2=t\gamma^{-1}(\lambda(b+c)+\mu c),f_3=t\gamma^{-1}b$ with $\mu,\lambda$ distinct roots of $x^2+\gamma x+1$.      \\ \hline
				$J_{3,31}^{\gamma\neq 0} \to J_{3,25}$        & $f_1=ta,f_2=tb,f_3=t^2c$\                                    \\ \hline
				$J_{3,31}^{\gamma\neq0} \to J_{3,26}$         & $f_1=(t+t\gamma^{-1})b+tc,f_2=\gamma^{-1}a+\gamma^{-1}c,f_3=tb$                          \\ \hline
				$J_{3,31}^{\gamma\neq0} \to   J_{3,27}$        &    $f_1=t^2b,f_2=t(b+c),f_3=a$\                    \\ \hline
				$J_{3,31}^{\gamma\neq 0}\to J_{3,28}^{\delta\neq 0}$ & $f_1=\beta(ra+c),\ f_2=ta,\ f_3=\beta tb$ with $\beta=\delta/\gamma$ and $r=1+\gamma/\delta$, where $\gamma\delta\neq0$.\ \\ \hline
				$J_{3,32} \to J_{3,25}$           & $f_1=ta,f_2=tb,f_3=t^2c$, 
				\\ \hline
				$J_{3,32} \to J_{3,26}$          & $f_1=tc,f_2=\alpha a+c,f_3=t\alpha b$, where $\alpha^2=\alpha+1$                                  \\ \hline
				$J_{3,32} \to J_{3,27}$  &  $f_1=t\alpha c,f_2=t\alpha b,f_3=a+\alpha b+\alpha c$, where $\alpha^2=\alpha+1$                   \\ \hline
				$J_{3,32} \to J^{\gamma\neq0}_{3,28}$        &   $f_1=\gamma\lambda^{-1} (a+\lambda b), f_2=tb, f_3=\gamma\lambda^{-1}tc$, $\lambda$ a root of $g(x)=(\gamma^2+1)x^2+\gamma^2 x+\gamma^2$  
				\\ \hline
				$J_{3,15}^1 \to J_{3,3}^1$        &  $f_1=tb,\ f_2=ta,\ f_3=t^2a^2$                          \\ \hline
				$J_{3,16}^{1,0} \to J_{3,3}^1$        &  $f_1=ta,f_2=t+b,f_3=t^2$                          \\ \hline
				$J_{3,16}^{1,0} \to J_{3,6}$        &  $f_1=1,f_2=ta,f_3=b$                          \\ \hline
				$J_{3,16}^{1,0} \to J_{3,26}$        &  $f_1=t+b,f_2=1+a,f_3=b$                          \\ \hline
				$J_{3,16}^{1,0} \to J_{3,15}^0$        &  $f_1=1+(1+t)a,f_2=t^2,f_3=b$                          \\ \hline
			\end{tabular}
		\end{adjustbox}
		\caption{Existence of deformations in $Jor_3$}\label{table2}
	\end{table}

	\begin{proposition}
		Regarding the Characteristic condition, we have the following conditions:
		\begin{itemize}
			\item For $x=s1+ua+vb\in J_{3,17}$, the corresponding characteristic polynomial is $\chi_{L_x}(\lambda)=\lambda^3+s\lambda^2+s^2\lambda+s^3$;
			\item For $x=s1+ua+vb\in J_{3,16}^{1,0}$, the corresponding characteristic polynomial is $\chi_{L_x}(\lambda)=\lambda^3+s\lambda^2+(s^2+u^2)\lambda+s^3+su^2$;
			\item For $x=sa+ub+vc\in J_{3,26}$, the corresponding characteristic polynomial is $\chi_{L_x}(\lambda)=\lambda^3+u\lambda^2$;
			\item For $x=sa+ub+vc\in J_{3,28}^\gamma$, the corresponding characteristic polynomial is $\chi_{L_x}(\lambda)=\lambda^3+\gamma s\lambda^2+s^2\lambda$;
		\end{itemize}
	\end{proposition}
	\begin{proof}
		In the order of the four items, the matrices of $L_x$ in the displayed bases are
		$$
		\begin{pmatrix}s&v&u\\u&s&0\\v&0&s\end{pmatrix},\quad
		\begin{pmatrix}s&u&0\\u&s&0\\v&0&s\end{pmatrix},\quad
		\begin{pmatrix}0&0&0\\0&0&0\\u&s+v&u\end{pmatrix},\quad
		\begin{pmatrix}0&0&0\\v&0&s\\u+\gamma v&s&\gamma s\end{pmatrix}.
		$$
		Expanding $\det(\lambda I-L_x)$ gives the four formulas. In particular, for $J_{3,16}^{1,0}$ it equals $(\lambda+s)((\lambda+s)^2+u^2)$, whose constant coefficient is $s^3+su^2$.
	\end{proof}
	
	\begin{theorem}
		For the algebraic variety $Jor_3$ over an algebraically closed field of characteristic two, the following holds:
		\begin{itemize}
			\item The irreducible components of $Jor_3$ are $\overline{(J_{3,13}^0)^{GL(V)}},\overline{(J_{3,16}^{1,1})^{GL(V)}},\overline{(J_{3,24})^{GL(V)}},\mathcal{Z}$, where
			$$\mathcal{Z}=\{A\in Jor_3 : x^2=0,\forall x\in A\}$$
			\item The rigid algebras are exactly $J_{3,13}^0,J_{3,16}^{1,1}$ and $J_{3,24}$.
		\end{itemize}
	\end{theorem}
	\begin{proof}
		Write $C_1=\overline{(J_{3,13}^0)^{GL(V)}}$,
		$C_2=\overline{(J_{3,16}^{1,1})^{GL(V)}}$ and
		$C_3=\overline{(J_{3,24})^{GL(V)}}$. The classification and the degenerations place every algebra with a nonzero square in the union $C_1\cup C_2\cup C_3$. Every remaining algebra belongs to $\mathcal Z$.
		Hence $\mathrm{Jor}_3=C_1\cup C_2\cup C_3\cup\mathcal Z$.
		The three orbit closures are irreducible, and their dimensions are
		$$
		\dim C_1=9,\qquad \dim C_2=8,\qquad \dim C_3=7.
		$$
		By Theorem \ref{zeropotentcomponent}, $\mathcal Z$ is irreducible of
		dimension $9$.
		
		No two of these four closed sets contain one another. Associativity
		holds throughout $C_1$, and the defining algebras of $C_2$ and $C_3$
		are not associative. The condition that every $L_{x^2}$ has the form $c(x)I$ holds
		on $C_2$, but fails in $J_{3,24}$: for its idempotent $a$, the operator
		$L_{a^2}=L_a$ has diagonal $(1,0,0)$. Using the dimension we conclude that the reverse containment is not possible. Each $C_i$ contains an algebra with a nonzero square,
		so $C_i\not\subseteq\mathcal Z$. Dimension gives
		$\mathcal Z\not\subseteq C_2,C_3$. Notice that $\mathcal Z\not\subseteq C_1$ because $J_{3,27}$ is not associative. Thus these are exactly the irreducible components.
		
		Each defining orbit is open in its orbit closure. It meets no other
		component, since such an intersection would force its entire orbit
		closure into that component. Its complement in the whole variety is
		therefore the union of its closed boundary and the other three
		components. Hence the three defining algebras are rigid.  There is no open orbit in $\mathcal Z$.  The classification of three-dimensional square-zero algebras shows that
		their isomorphism classes are represented by
		$$
		J_{3,1},\ J_{3,25},\ J_{3,26},\ J_{3,27},\
		J_{3,28}^{\gamma\neq0},\ J_{3,28}^{0},\
		J_{3,29},\ J_{3,30},\
		J_{3,31}^{\gamma\neq0},\ J_{3,31}^{0},\ J_{3,32}.
		$$
		Their automorphism-group dimensions are, respectively,
		$$
		9,\ 6,\ 4,\ 6,\ 4,\ 4,\ 1,\ 2,\ 1,\ 3,\ 1.
		$$
		Thus every orbit in $\mathcal Z$ has dimension
		$$
		\dim\bigl(\mathrm{GL}_3(K)\cdot A\bigr)
		=9-\dim\operatorname{Aut}(A)\leq 8.
		$$
		Since $\mathcal Z$ is irreducible of dimension $9$, an open orbit in
		$\mathcal Z$ would have dimension $9$. Therefore $\mathcal Z$ has no
		open orbit.
		Finally, every other algebra with a nonzero square
		has a proper incoming degeneration from one of the three defining
		orbits. Such an algebra cannot have an open orbit. This proves rigidity.
	\end{proof}
	
	\section*{Acknowledgements}
	The first author was financed by Secihti-80325 PY 1449. The fourth author was financed in part by the Coordenação de Aperfeiçoamento de Pessoal de Nível Superior - Brasil (CAPES) - Finance Code 001.

\end{document}